\documentclass[12pt]{amsart}
\usepackage{amssymb}
\usepackage{graphicx}

\usepackage{mathrsfs}

\def\e{\text{\rm e}}

\newtheorem{theorem}{Theorem}
\newtheorem{lemma}{Lemma}

\newtheorem{prop}{Proposition}

\theoremstyle{definition}
\newtheorem{df}{Definition}

\newtheorem{remark}{Remark}

\title{Optimal Consumption Problem in a Diffusion Short-Rate Model}

\begin{document}

\maketitle

\begin{center}
Daniel Synowiec\footnote{Faculty of Applied Mathematics,
        AGH University of Science and Technology, al. Mickiewicza 30,
        30-059 Cracow, Poland ({\tt synowiec@agh.edu.pl}).}
\end{center}

\begin{abstract}
We consider a problem of an optimal consumption strategy on the
infinite time horizon when the short-rate is a diffusion process.
General existence and uniqueness theorem is illustrated by the
Vasicek and so-called invariant interval models. We show also that
when the short-rate dynamics is given by a Brownian motion or a
geometric Brownian motion, then the value function is infinite.
\end{abstract}

\begin{center}
\end{center}

\begin{center}
\end{center}

\pagestyle{myheadings} \thispagestyle{plain} \markboth{D.
SYNOWIEC}{OPTIMAL CONSUMPTION IN A SHORT-RATE MODEL}

\section{Introduction}
Let $r_t$ be the short-rate (i.e. the rate offered by a bank) at
time $t \geq 0$. Assume that $(r_t)$ satisfies the following
stochastic differential equation
\begin{equation} \label{r}
\begin{aligned}
d r_t &= \mu(r_t) dt + \sigma(r_t) dW_t,\\
r_0 &= r,
\end{aligned}
\end{equation}
where $(W_t)$ is a one-dimensional Brownian motion defined on a
filtered probability space
$(\Omega,\mathcal{F},(\mathcal{F}_t),\mathbb{P})$.

Let us denote by  $V^{(C;r,v)}_t$  the capital at time $t$ of a
bank account owner whose consumption rate is $C$ and whose wealth
at time 0 is $v > 0$. Then
$$
d V^{(C;r,v)}_t =\left( r_t V^{(C;r,v)}_t-C_t \right)dt ,\qquad
V^{(C;r,v)}_0 = v.
$$
Let
\begin{equation} \label{tau}
\tau_A^{(C;r,v)}=\inf\left\{t\ge 0: V^{(C;r,v)}_t = 0\right\}
\end{equation}
be the bankruptcy time.

In the paper it is assumed that any consumption rate $C$ is
progressively measurable and non-negative. The space of all
consumption rates is denoted by  $\mathcal{U}$.

Given a discount factor $\gamma \ge 0$ and an exponent $\alpha \in
(0,1)$ of the power utility function, we are concerned with the
following problems:

\bigskip
\noindent \textbf{Problem A} Given $r$ and $v$, find a consumption
rate $\hat{C}^{(r,v)}\in \mathcal{U}$ which maximizes the
performance functional, that is,
$$
J(\hat{C}^{(r,v)};r,v)=\sup_{C\in \mathcal{U}} J(C;r,v),
$$
where
\begin{equation} \label{J}
J(C;r,v):= \mathbb{E}^r \int_0^{\tau_A} \e ^{-\gamma t} C^\alpha_t
dt,
\end{equation}
$\tau_A = \tau_A^{(C;r,v)}$, and $\mathbb{E}^r$ is the conditional
expectation $\mathbb{E}\left( \cdot \, \vert r_0=r\right)$. A
solution to this problem is  given in Proposition \ref{prop-A} and
Theorem \ref{t-sol-a} below.

\bigskip
\noindent \textbf{Problem B} It is reasonable to assume that one
keeps his money in the bank account as long as the interest rate
$r_t$ is positive. Under this assumption the performance
functional is given by
$$
J_B(C;r,v):= \mathbb{E}^r \left[ \int_0^{\tau_B} \e ^{-\gamma t}
C^\alpha_t dt + \e^{-\gamma \tau_B} V_{\tau_B}^\alpha \cdot
\chi_{\{\tau_B < \infty\}} \right],
$$
where $\tau_B = \tau_A^{(C;r,v)} \wedge \tau_0^r$ , $\tau_0^r =
\inf\left\{t\ge 0: r_t = 0\right\}$ and $V = V^{(C;r,v)}$.
Clearly, if $r \leq 0$, then $\tau_B =0$.  The goal is now to find
a consumption rate which maximizes  $J_B$ (see Proposition
\ref{prop-B} and Theorem \ref{t-sol-b}).

\bigskip
\noindent \textbf{Problem C} Let $p(t,\theta)$ be the price at
time $t$ of a zero-coupon bond that pays off $1$ at time $\theta$.
If one may also invest in zero-coupon bonds then the wealth
dynamics is formally given by
\begin{equation} \label{V-bond}
\begin{aligned}
d V_t^{(u;r,v)} &= \left( \eta_t r_t V_t^{(u;r,v)}-C_t \right)dt\\
&\quad + (1-\eta_t) V_t^{(u;r,v)} \int_0^\infty
\frac{dp(t,\theta)}{p(t,\theta)}
\psi(t,\theta) d\theta,\\
V_0^{(u;r,v)} &= v,
\end{aligned}
\end{equation}
where $u = (C,\eta,\psi)$, $\psi$ is  the density of the
distribution of investments in bonds with various terminal time
$\theta$. The aim is to maximize the performance functional $J$
given by \eqref{J} with $\tau_A = \tau_A^{(u;r,v)}$.

\bigskip Problems A, B and C defined above are particular cases of
an investor problem, various types of which has been investigating
since 1970's (see \cite{Merton} and \cite{Merton-art}). However,
most of them are concerned with investment in a bank account
(usually on a constant rate) and a finite number of stocks. If one
can invest in a bank account and zero-coupon bonds, then the
investor problem is more difficult to solve. The reason is that
there can be an infinite number of bonds, since the time of
maturity $\theta$ can take an infinity of values. Furthermore, the
set of admissible strategies does not contain "buy and hold"
strategy, i.e. one must convert bond to cash at maturity.

The type of an investor who can invest his money in bonds has been
recently studied in \cite{Ekeland-Taflin}, \cite{Korn-Kraft} and
\cite{Ringer-Tehranchi}. Contrary to our paper, authors of
\cite{Ekeland-Taflin}, \cite{Korn-Kraft} and
\cite{Ringer-Tehranchi} examined the portfolio problem without
possibility of consumption and with a finite time horizon. On the
other hand, in \cite{Ekeland-Taflin} and \cite{Ringer-Tehranchi}
it is assumed that the dynamics of the instantaneous forward rate
is given and that the performance function is defined under a real
measure. More references can be find in the survey paper
\cite{Zariphopoulou}.

In the paper we use the Hamilton--Jacobi--Bellman approach,
whereas in \cite{Ekeland-Taflin} and \cite{Ringer-Tehranchi}
convex duality is used.

\section{Preliminaries}

In the paper, it is assumed that \eqref{r} defines a Markov family
on an open subinterval $\mathcal {O} \subseteq \mathbb{R}$, which,
in particular, means that $\mathcal{O}$ is invariant for
\eqref{r}; that is, $r_0\in \mathcal{O}$ implies that $r_t\in
\mathcal{O}$ for all $t\ge 0$. Moreover, it is assumed that
$\sigma \in C^2(\mathcal {O})$, $\mu \in C^1(\mathcal {O})$, their
first derivatives are bounded on $\mathcal{O}$, and that the
diffusion is non-degenerate, i.e. $\sigma(r) \neq 0$ for all $r
\in \mathcal{O}$.

The \textit{value function} for one of the listed above problems
is the maximum of the corresponding performance functional over
the set of admissible controls. We will show that the value
functions are very regular, namely $C^2$ in $r$. Let
\begin{equation}\label{generator}
Q f(r) := \frac 12 \sigma^2(r) f''(r)  + \mu(r) f'(r),
\end{equation}
be the \textit{formal generator} of the diffusion given by \eqref{r}.

The results below have the form of the verification theorem for
stochastic control problems. For similar results see e.g.
\cite{Fleming-Soner}, \cite{Oksendal} or \cite{Oksendal-Sulem}.

\begin{prop}\label{prop-A}
Let $K \in C^2(\mathcal{O})$ be such that
\begin{equation}\label{K-hjb}
Q K(r) + (\alpha r-\gamma) K(r) + (1-\alpha)
K^\frac{\alpha}{\alpha-1}(r)=0,
\end{equation}
for every $r \in \mathcal{O}$. Then $\Phi(r,v)=K(r)v^\alpha$ is
the value function for Problem A, whenever for any $C \in
\mathcal{U}$ and $r \in \mathcal{O}$,
\begin{equation}\label{Condition}
\lim_{n \to \infty} \mathbb{E}^r \e^{-\gamma \tau_n}
\Phi(r_{\tau_n},V_{\tau_n}) =0,
\end{equation}
where $(r_t)$ is given by \eqref{r}, $\tau_n = n \wedge
\tau_A^{(C;r,v)}$ and $V = V^{(C;r,v)}$. The optimal consumption
is given in the feedback form
\begin{equation}\label{C-opt}
\hat{C} = K^{\frac{1}{\alpha -1}}v.
\end{equation}
\end{prop}

\noindent {\bf Proof} Taking into account the dynamics of $(V_t)$
and the form of performance functional we see that
$$
\Phi(r,v) = K(r) v^\alpha
$$
for a certain function $K$. The Hamilton--Jacobi--Bellman equation
(see e.g. \cite{Fleming-Soner}, \cite{Oksendal-Sulem}) for $K$ is
$$
\sup_{C \geq 0}\left\{ -\gamma K(r) v^\alpha + Q K(r) v^\alpha +
\alpha (rv-C) v^{\alpha-1} K(r) + C^\alpha\right\} =0.
$$
The supremum is attained at $\hat{C}$ given by \eqref{C-opt}.
Hence, $K$ satisfies \eqref{K-hjb} and the HJB verification
theorem (see \cite{Oksendal}, \cite{Oksendal-Sulem}) gives us the
claim. $_\square$

\bigskip In Problem B we have to assume that $0 \in \mathcal{O}$. If not, Problem B can
be reduced to Problem A. Let $\mathcal{O}^+ = \mathcal{O} \cap [0,
\infty)$ and $\mathcal{O}^{++} = \mathcal{O} \cap (0, \infty)$.
With a similar proof as above we have the following proposition
concerning Problem B.

\begin{prop}\label{prop-B}
Let $K \in C^2(\mathcal{O}^{++}) \cap C(\mathcal{O}^+)$ satisfy
\eqref{K-hjb}. Then
$$
\Phi(r,v)= \left\{
\begin{array}{ll}
    v^\alpha, & \hbox{$r \in \mathcal{O} \setminus \mathcal{O}^{++}$}, \\
    K(r)v^\alpha, & \hbox{$r \in \mathcal{O}^{++}$}, \\
\end{array}
\right.
$$
is the value function for Problem B, whenever for any $C \in
\mathcal{U}$ and $r \in \mathcal{O}^{++}$,
\begin{equation}\label{Condition-B}
\lim_{n \to \infty} \mathbb{E}^r \e^{-\gamma \tau_n}
\Phi(r_{\tau_n},V_{\tau_n}) = \mathbb{E}^r \e^{-\gamma \tau_B}
V_{\tau_B}^\alpha \chi_{\{\tau_B<\infty\}},
\end{equation}
where $\tau_B = \tau_B^{(C;r,v)}$, $\tau_n^{(C;r,v)} = n \wedge
\tau_B$ and $V = V^{(C;r,v)}$. The optimal consumption is given in
the feedback form \eqref{C-opt}.
\end{prop}

\bigskip Note that $K$ satisfies a non-linear, non-Lipschitz second
order differential equation, but $K$ is not defined as a solution
to the Cauchy problem. The goal of the paper is to prove the
existence of the solution satisfying appropriate boundary
conditions and to find an approximating scheme for $K$.

\section{Solution to Problem C}

In Problem C we assume that  $\mathbb{P}$ is a martingale measure.
Then
$$
p(t,\theta) = \mathbb{E}\left( \e^{-\int_t^\theta r_s
ds}|\mathcal{F}_t \right).
$$
Since $(r_t)$ is a Markov process, $p(t,\theta)=
\nu^\theta(t,r_t)$ is a function of $t,\theta$ and $r_t$. Thus we
can rewrite \eqref{V-bond} as follows
\begin{equation}\label{V-bond-mart}
dV_t = (r_t V_t-C_t)dt + (1-\eta_t) V_t \sigma(r_t) \int_0^\infty
\frac{\frac{\partial}{\partial r}
\nu^\theta(t,r_t)}{\nu^\theta(t,r_t)} \psi(t,\theta) d\theta dW_t
\end{equation}
with $V = V^{(u;r,v)}$, whenever $\nu^\theta$ is differentiable
with respect to $r$.

In Problem C the performance function is given by \eqref{J}, the
class of admissible controls $\mathcal{U}$ consists of tuples
$(C,\eta,\psi(\cdot,\theta)_{\theta\ge 0})$ of progressively
measurable processes, such that $C_t$ is non-negative,
\eqref{V-bond-mart} is well defined and
$$
\int_0^\infty \psi(t,\theta) d\theta =1, \qquad  \qquad
\psi(t,\theta) \equiv 0, \quad \forall \theta \leq t.
$$
Note that neither $\eta$ nor $\psi$ have to be non-negative.

Define
\begin{equation}\label{beta}
\Upsilon_s = \int_0^\infty \frac{\frac{\partial}{\partial r}
\nu^\theta(s,r_s)}{\nu^\theta(s,r_s)} \psi(s,\theta)  d\theta
\quad\text{and}\quad \beta_s = (1-\eta_s) \Upsilon_s.
\end{equation}
Then we can rewrite \eqref{V-bond-mart} in the form
\begin{equation}\label{V-bond-mart-beta}
dV_t = (r_t V_t-C_t)dt + \beta_t V_t \sigma(r_t) dW_t.
\end{equation}
Since we assumed that we are given dynamics of $(r_t)$, and the
performance functional is under a martingale measure, we can treat
the investor portfolio as the one consisted of the bank account
and one other instrument with price $S_t$ given by
\begin{equation}\label{dS}
dS_t = S_t \left( r_t dt + \sigma(r_t) \Upsilon_t dW_t \right),
\qquad S_0=1.
\end{equation}
Note that given \eqref{V-bond-mart-beta} and \eqref{dS}, we have
to assume that for any $t>0$,
$$
\int_0^t \Upsilon_s^2 ds < \infty, \qquad \int_0^t \beta_s^2 ds <
\infty,\qquad \mathbb{P}-\text{a.s.}
$$
Therefore the dynamics of the wealth of the investor is given by
$$
d V_t = (\eta_t r_t V_t - C_t)dt + (1-\eta_t) V_t dS_t/S_t,
$$
which is equivalent to \eqref{V-bond} and
\eqref{V-bond-mart-beta}. Thus the number of instruments is finite
and the same approach as in \cite{Karatzas-Shreve} can be taken.

Theorem \ref{t-sol-c} below, giving a solution to Problem C, was
formulated and proven in \cite{Karatzas-Shreve} under much weaken
conditions. Here we present another proof, based on Proposition
\ref{prop-C} below. We restrict our attention to the value
function of the problem. We refer the reader to
\cite{Karatzas-Shreve} for details on the optimal control
(portfolio).

\begin{prop}\label{prop-C}
If $K \in C^2(\mathcal{O})$ satisfies
\begin{equation}\label{K-hjb-bonds}
Q K(r) + (\alpha r-\gamma) K(r) + (1-\alpha)
K^\frac{\alpha}{\alpha-1}(r) + \frac{\alpha
\sigma^2(r)}{2(1-\alpha)}
 \frac{(K'(r))^2}{K(r)}=0,
\end{equation}
then $\Phi(r,v)=K(r)v^\alpha$ is the value function for Problem C,
whenever for any $u \in \mathcal{U}$ and $r \in \mathcal{O}$,
\begin{equation}\label{Condition-C}
\lim_{n \to \infty} \mathbb{E}^r \e^{-\gamma \tau_n}
\Phi(r_{\tau_n},V_{\tau_n}) = 0,
\end{equation}
where $\tau_n = n \wedge \tau_A^{(u;r,v)}$ and $V = V^{(u;r,v)}$.
The optimal consumption $\hat C$ and the optimal factor $\hat
\beta$ defined in \eqref{beta} are given by
\begin{equation}\label{C-b-opt}
\hat{C} = K^{\frac{1}{\alpha -1}}v , \qquad \hat{\beta} =
\frac{K'}{(1-\alpha)K}.
\end{equation}
\end{prop}
\noindent {\bf Proof} Again $\Phi(r,v) = K(r) v^\alpha$ for a
certain function $K$. The HJB equation for $K$ is
\begin{equation}\label{hjb-beta}
\begin{aligned}
Q K(r) v^\alpha + (\alpha r-\gamma) K(r) v^\alpha + \sup_{C \geq
0}\left\{C^\alpha - C \alpha K(r) v^{\alpha-1} \right\}&\\
+\sup_{\beta}\left\{\alpha \beta \sigma^2(r) K'(r) v^\alpha +
\frac{\alpha(\alpha-1)}{2} \beta^2 \sigma^2(r) K(r) v^\alpha
\right\} &=0,
\end{aligned}
\end{equation}
and the claim  follows from the HJB verification theorem.
$_\square$

\bigskip Since the value function  $\Phi(r,v)$ is a
non-decreasing positive function of both arguments $r$ and $v>0$,
we see that $K$ is non-decreasing and  positive. Then the optimal
$\hat{\beta}$ in \eqref{C-b-opt} is positive. Note that $\eta \leq
1$ and $\psi \geq 0$ whilst the short-selling is forbidden. Then
the condition
\begin{equation}\label{vega}
\frac{\partial}{\partial r} \nu^\theta(t,r_t) \leq 0,
\end{equation}
which holds e.g. in Vasicek and CIR models, implies that if the
short-selling is forbidden then necessarily $\Upsilon_t$ and
$\beta_t$ given by \eqref{beta} are non-positive. Thus the
supremum over $\beta \leq 0$ in equation \eqref{hjb-beta} is
attained at 0. Hence, if the short-selling is forbidden and
\eqref{vega} holds, then Problem C reduces to Problem A.

It is worth mentioning that given $\hat{\beta}$ we do not have
unambiguous solution to Problem C, i.e. we do not obtain
unambiguous pair $(\eta,\psi)$. However we may choose arbitrary
$\psi$ such that $(C,\eta,\psi) \in \mathcal{U}$ and then we
derive an optimal $\hat{\eta}$ from \eqref{beta}. For example we
may set $\psi(t,\theta)=\varsigma \e^{-\varsigma (\theta -t)}
\cdot \chi_{\{ t < \theta \}}$ for some $\varsigma >0$.

\bigskip The following result will be used to show the regularity of
the value function. Let
\begin{equation} \label{N}
N(r) := \mathbb{E}^r \int_0^\infty \e^{\frac{1}{1-\alpha}(-\gamma t
+ \alpha \int_0^t r_s ds)}dt, \qquad r\in \mathcal{O}.
\end{equation}

\begin{prop}\label{prop-NinC2}
If $N(r) < \infty$ for every $r \in \mathcal{O}$ and
\begin{equation} \label{N2}
{\bf E}^r \int_0^\infty \e^{\frac{2}{1-\alpha}(-\gamma t + \alpha
\int_0^t r_s ds)}dt < \infty, \qquad \forall r\in \mathcal{O},
\end{equation}
then $N \in C^2(\mathcal{O})$ and
\begin{equation} \label{N-hjb}
Q N(r) + \frac{\alpha r-\gamma}{1-\alpha} N(r) + 1 =0, \qquad r\in
\mathcal{O}.
\end{equation}
\end{prop}

\noindent {\bf Proof} Let $\mathcal{O} = (a,b)$, where $a \geq
-\infty$ and $b \leq \infty$. Then, by (\ref{N2}),
$$
N_n(r) = {\bf E}^r \int_0^{\tau_{\pm n}^r}
\e^{\frac{1}{1-\alpha}(-\gamma t + \alpha \int_0^t r_s ds)}dt
$$
is a solution to the boundary problem
$$
\left\{
\begin{array}{ll}
    Q N_n(r) + \frac{\alpha r - \gamma}{1-\alpha} N_n(r) = -1,\\
    N_n(a_n) = N_n(b_n) = 0,\\
\end{array}
\right.
$$
where
$$
a_n = \left\{
\begin{array}{ll}
    a + 1/n, & \hbox{$a > -\infty$}, \\
    -n, & \hbox{$a = -\infty$}, \\
\end{array}
\right. \qquad b_n = \left\{
\begin{array}{ll}
    b - 1/n, & \hbox{$b < \infty$}, \\
    n, & \hbox{$b = \infty$}, \\
\end{array}
\right.
$$
and $\tau_{\pm n}^r = \inf\{ t \geq 0 \colon r_t \notin [a_n,b_n]
\}$. Since we assumed that $\mathcal{O}$ is invariant for
\eqref{r}, then $\lim_{n \to \infty} \tau_{\pm n}^r = \infty$ for
any $r \in \mathcal{O}$ and consequently $N(r) = \lim_{n \to
\infty} N_n(r)$. Thus $N$ is a weak solution (see Definition
\ref{def-weak} in Section \ref{sec-tools}) to \eqref{N-hjb}, and
by Lemma \ref{l-AH}, $N \in C^2(\mathcal{O})$. Hence, it is a
strong solution to \eqref{N-hjb}. $_\square$

The result below says that the function $K$ appearing in the
identity $\Phi(r,v) =K(r)v^\alpha$  for  the value function equals
$N^{1-\alpha}$.

\begin{theorem}\label{t-sol-c}
Let assumptions of Proposition \ref{prop-NinC2} hold. Assume
additionally that for any $u \in \mathcal{U}$ and $r \in
\mathcal{O}$,
\begin{equation}\label{Condition-N}
\lim_{n \to \infty} \mathbb{E}^r \e^{-\gamma \tau_n}
N^{1-\alpha}(r_{\tau_n}) V_{\tau_n}^\alpha = 0,
\end{equation}
where $\tau_n = n \wedge \tau_A^{(u;r,v)}$ and $V = V^{(u;r,v)}$.
Then $\Phi(r,v) = N^{1-\alpha}(r) v^{\alpha}$ is the value
function for Problem C.
\end{theorem}
\noindent {\bf Proof} By elementary calculus, \eqref{N-hjb} is
equivalent to \eqref{K-hjb-bonds} for $K(r) = N^{1-\alpha}(r)$.
Condition \eqref{Condition-N} implies \eqref{Condition-C} and we
conclude by Proposition \ref{prop-C}. $_\square$

\section{Solution to Problem A}

This section contains one of the main result of the paper. It
provides the existence and approximating scheme for the solution
$K$ to the HJB equation \eqref{K-hjb} for Problem A.  In its
formulation  $(E, \|\cdot\|_E)$ is a Banach space of continuous
functions on $\mathcal{O}$.

We will need the following hypotheses:

\smallskip \noindent \textbf{(H.1)} For any fixed $t \geq 0$, $r \in \mathcal{O}$, $\varphi \in E$ and
any sequence $\{ T_n \}$ of stopping times, the sequences of
random variables
$$
\left\{ \varphi(r_{t \wedge T_n}) \e^{ \alpha \int_0^{t \wedge
T_n} r_s ds}  \right\}_{n \in \mathbb{N}} \quad {\rm and} \quad
\left\{ \int_0^{t \wedge T_n} \varphi(r_{s}) \e^{ \alpha
\int_0^{s} r_u du} ds \right\}_{n \in \mathbb{N}}
$$
are uniformly integrable with respect to $\mathbb{P}^r =
\mathbb{P}\left( \cdot \, \vert r_0=r\right)$.

\smallskip \noindent \textbf{(H.2)} For any Lipschitz continuous
bounded function $f \colon [0,\infty) \mapsto [0,\infty)$ and for
any non-negative $\phi \in E$, one has $f(\phi)\in E$.

\smallskip
\noindent\textbf{(H.3)} The family  $(P_t, t\geq 0)$  of linear
operators
\begin{equation}\label{P}
P_t \varphi(r) = \mathbb{E}^r \varphi(r_t) \e^{\alpha \int_0^t r_s
ds}, \qquad r \in \mathcal{O},
\end{equation}
forms  a $C_0$-semigroup on $E$.

\begin{remark}
In examples $\mathcal{O}=\mathbb{R}$ and
$$
E = \{\varphi \in C(\mathbb{R}) \colon \lim_{|r| \to \infty} |\varphi(r)| \e^{-\frac\alpha b |r|} = 0 \}
$$
or $\mathcal{O}$ is a bounded interval and $E$ is the space
$UC(\mathcal{O})$ of uniformly continuous functions on $\mathcal
{O}$. Moreover, we will show in Lemma \ref{l-Generator} that the
generator $(\mathcal{A}, \mathcal{D}(\mathcal {A}))$ of $(P_t)$ is
given by
$$
\mathcal{D}(\mathcal{A}) = \{ \varphi \in C^2(\mathcal{O}) \cap E \colon  A \varphi \in E\},
$$
and  $\mathcal{A} \varphi = A \varphi$  for all $\varphi \in \mathcal{D}(\mathcal{A})$, where  $A$ is the  differential operator
\begin{equation} \label{A}
 A\varphi(r) = Q\varphi(r)+\alpha r \varphi(r).
 \end{equation}
\end{remark}

\smallskip
We note that condition \textbf{(H.1)} will be needed only in the
proof of the inclusion $\{ \varphi \in C^2(\mathcal{O}) \cap E
\colon  A \varphi \in E\} \subset \mathcal{D}(\mathcal{A})$.

\smallskip
Recall that $N$ is a function defined by \eqref{N}. The following
hypothesis will be needed in the proof  that $N^{1-\alpha}$ is a
supersolution to the HJB equation \eqref{K-hjb}, such that
$N^{1-\alpha} \in {\mathcal D}({\mathcal A})$ and $\Phi_N(r,v) =
N^{1-\alpha}(r) v^\alpha$ satisfies the boundary condition
(\ref{Condition}). For more details see Definition \ref{Def} and
Remark \ref{RSuper}.

\smallskip \noindent \textbf{(H.4)}
For any $r \in \mathcal{O}$,
\begin{equation} \label{boundary}
\lim_{t \to \infty} \mathbb{E}^r \e^{-\gamma t + \alpha \int_0^t
r_s ds} N^{1-\alpha}(r_t) = 0
\end{equation}
and for any stopping time $\tau_A^{(C;r,v)}$,
\begin{equation} \label{bnd-uni-int}
\left\{ \e^{-\gamma\tau_n + \alpha \int_0^{\tau_n} r_s ds}
N^{1-\alpha}(r_{\tau_n}) \chi_{\{\tau_A^{(C;r,v)} < \infty\}}
\right\}_{n \in \mathbb{N}}
\end{equation}
is uniformly integrable, where $\tau_n = n \wedge
\tau_A^{(C;r,v)}$.

Moreover, $N^{1-\alpha}\in C^2(\mathcal {O})\cap E$ and $AN^{1-\alpha}
\in E$, where $A$ is defined by \eqref{A}.

\smallskip
For any $m>0$, define
\begin{equation}\label{Fm}
F_m(x) = \left\{
\begin{array}{ll}
    (1-\alpha) x^{\frac{\alpha}{\alpha-1}}, & \hbox{$x > m^{\alpha-1}$}, \\
    m^\alpha - \alpha m x, & \hbox{$0 \leq x \leq m^{\alpha-1}$}. \\
\end{array}
\right.
\end{equation}

Recall that $\mathcal{A}$ is the generator of the semigroup
$(P_t)$. We denote by $\varrho(\mathcal{A}-\gamma)$ the resolvent
set of $\mathcal{A}-\gamma$. The proof of the following theorem is
postponed to Section \ref{sec-proofA}.

\begin{theorem}\label{t-sol-a}
Assume that \textbf{(H.1)}, \textbf{(H.2)}, \textbf{(H.3)} and
\textbf{(H.4)} are fulfilled. Then there is a solution $K$ to
\eqref{K-hjb} with condition \eqref{Condition}. Moreover, $K(r)
\leq N^{1-\alpha}(r)$, $r \in \mathcal{O}$. Finally,  for any
sequence $\{\lambda_m\} \subset \varrho(\mathcal{A}-\gamma)$ such
that for any $m
> 0$,
\begin{equation}\label{F+l}
the \; mapping \; \; [0,\infty) \ni x \mapsto F_m(x) + \lambda_m x
\; \; is \; non {\emph-} decreasing,
\end{equation}
one has
$$
K(r) = \lim_{m \to \infty} \lim_{n \to \infty} K^m_n(r), \qquad r
\in \mathcal{O},
$$
where $\{K^m_n\}$ is a non-decreasing sequence of both $m$ and
$n$, defined as follows
$$
\begin{aligned}
& K^m_0 = 0,\\
& K^m_{n+1} = (\lambda_m +\gamma -\mathcal{A})^{-1}(F_m(K^m_n) +
\lambda_m K^m_n).
\end{aligned}
$$
\end{theorem}

\begin{remark}\label{rem-Fm}
Since $F_m$ are Lipschitz continuous, then the function $x \mapsto
F_m(x) + \lambda x$ is non-decreasing for $\lambda$ large enough.
Thus there is a sequence $\{ \lambda_m \} \subset
\varrho(\mathcal{A}-\gamma)$ such that the functions $x \mapsto
F_m(x) + \lambda_m x$ are non-decreasing. Furthermore, from
$C_0$-semigroup property of $(P_t)$ guaranteed by \textbf{(H.3)}
we get
$$
\| P_t \varphi \|_E \leq M \e^{\vartheta t} \| \varphi \|_E
$$
for some $\vartheta$ and $M>0$. Then $(\vartheta, \infty) \subset
\varrho(\mathcal{A})$ and setting any $\varepsilon_1 >0$ and
$\varepsilon_2 \geq 0$ we may define
\begin{equation}\label{lambda}
\lambda_m = \max \{ \vartheta - \gamma + \varepsilon_1, \alpha m +
\varepsilon_2 \}.
\end{equation}
\end{remark}

\begin{remark}
We will show in Sections \ref{SVasicek} and \ref{exa-bnd}, that
the assumptions of the Theorem \ref{t-sol-a} are satisfied if
$(r_t)$ is an Ornstein--Uhlenbeck process (the so-called Vasicek
model) or $\mathcal{O}$ is bounded. We will show in Section
\ref{SInfinite} that if $(r_t)$ is either a Brownian motion or a
geometric Brownian motion then the value function for Problem $A$
is infinite.
\end{remark}

\section{Analytical tools}\label{sec-tools}

This section provides some useful analytical tools. Let us
consider a second order differential operator
$$
Du(x) = a_2(x) u''(x) + a_1(x) u'(x) + a_0(x) u(x)
$$
with $a_i \in C^i(\mathcal{O})$ and $a_2 \neq 0$ in $\mathcal{O}$. We denote by
$$
D^* u(x) =  (a_2(x)u(x))'' - (a_1(x) u(x))' + a_0(x) u(x)
$$
the formally adjoint operator. We denote by $L^1_{loc}(\mathcal {O})$ the space of all locally integrable  functions on $\mathcal O$.

\begin{df}\label{def-weak}
Let $f,u\in L^1_{loc}(\mathcal {O})$. We call $u$ a \textit{weak
solution} to the equation $Du=f$ if
$$
\int_{\mathcal{O}} u(x) D^* \varphi(x) dx = \int_{\mathcal{O}}
f(x)\varphi(x) dx,\qquad \forall\, \varphi \in C_0^\infty(\mathcal{O}).
$$
\end{df}

Let $\mathcal{G}$ be an open subset of $\mathbb{R}$. The following
result holds only in dimension $1$. For a counterexample in case
of ${\mathcal O} \subseteq \mathbb{R}^2$ see
\cite{Gilbarg-Trudinger}.

\begin{lemma}\label{l-AH}
Assume that $H \colon \mathcal{G} \mapsto \mathbb{R}$ is a
continuous function and $u \in L^1_{loc}(\mathcal {O})$  such that
$u({\mathcal O}) \subseteq \mathcal{G}$, is a weak solution to
\begin{equation}\label{AH}
Du = H(u).
\end{equation}
Then $u \in C^2(\mathcal{O})$, i.e. $u$ is a strong solution to \eqref{AH}.
\end{lemma}
\noindent \textbf{Proof} We may rewrite \eqref{AH} in the form
\begin{equation}\label{aiH}
(a_2 u' + (a_1-a_2')u)' = H(u) - (a_2'' - a_1' +a_0)u,
\end{equation}
where we skip argument $x$ and all derivatives of $u$ are in the
weak sense. We can use the following fact. Assume that $\xi$ is a
distribution whose derivative is a function $h \in
L^1_{loc}(\mathcal {O})$. Then $\xi$ is a function and
$$
\xi(x) = \zeta + \int_\triangle^x h(y) dy,
$$
for some finite $\triangle \in \mathcal{O}$ and $\zeta \in \mathbb{R}$.
Applying this observation to \eqref{aiH} we obtain
$$
a_2 u' + (a_1-a_2')u = \zeta_1 + \int_\triangle^x (H(u) - (a_2'' -
a_1' +a_0)u) dy,
$$
where the r.h.s. is continuous, since integrand is locally
integrable. Thus
$$
u' = \frac{\zeta_1}{a_2} + \frac{\int_\triangle^x (H(u) - (a_2'' -
a_1' +a_0)u) dy}{a_2} - \frac{(a_1-a_2')u}{a_2}
$$
and $u' \in L^1_{loc}(\mathcal {O})$. Using the same argument
again we have
$$
u(r) = \zeta_2 + \int_\triangle^r \left( \frac{\zeta_1}{a_2} +
\frac{\int_\triangle^x (H(u) - (a_2'' - a_1' +a_0)u) dy}{a_2} -
\frac{(a_1-a_2')u}{a_2} \right) dx
$$
and $u \in C$, since integrand is locally integrable. Having shown
that $u \in C(\mathcal{O})$, we see that the integrand is
continuous, which implies $u \in C^1(\mathcal{O})$. Now we
conclude that integrand is of class $C^1$ and consequently $u \in
C^2(\mathcal{O})$. $_\square$

\smallskip
Recall that $A$ is a differential operator given by \eqref{A}. We
denote by $(\mathcal{A}, \mathcal{D}(\mathcal {A}))$ the generator
of the $C_0$-semigroup $(P_t)$ defined  by \eqref{P} on the Banach
space $E$, see \textbf{(H.3)}.

\begin{lemma}\label{l-Generator}
We have
$$
\mathcal{D}(\mathcal{A}) = \{ \varphi \in C^2(\mathcal{O}) \cap E \colon  A \varphi \in E\},
$$
and  $\mathcal{A} \varphi = A \varphi$  for all $\varphi \in \mathcal{D}(\mathcal{A})$.
\end{lemma}
\noindent \textbf{Proof} Write $\mathcal{E} = \{ \varphi \in
C^2(\mathcal{O}) \cap E \colon  A \varphi \in E\}$.

\smallskip \textbf{Step 1.} Here we will show that $\mathcal{D}(\mathcal{A}) \subset
\mathcal{E}$. Let $\varphi \in \mathcal{D}(\mathcal{A})$. First we
will show that  $\langle \mathcal{A} \varphi, \psi \rangle =
\langle \varphi, A^* \psi \rangle$ for any $\psi \in
C_0^\infty(\mathcal{O})$. We have
$$
\begin{aligned}
\langle \mathcal{A} \varphi, \psi \rangle &= \lim_{t \downarrow 0}
\frac1t \int_\mathcal{O} (P_t \varphi(x) - \varphi(x)) \psi(x)
dx\\
&=  \lim_{t \downarrow 0} \frac1t \int_\mathcal{O}
\int_\mathcal{O} p_t(x,y) \varphi(y) \psi(x) dx dy - \lim_{t
\downarrow 0} \frac1t \int_\mathcal{O}  \varphi(y) \psi(y) dy,
\end{aligned}
$$
where $p_t(x,y)$ is a transition density function of process
$(r_t)$, which exists due to the non-degeneration of the diffusion
coefficient. Hence
$$
\begin{aligned}
\langle \mathcal{A} \varphi, \psi \rangle &= \lim_{t \downarrow 0}
\frac1t \int_\mathcal{O} \left( \int_\mathcal{O} p_t(x,y) \psi(x)
dx - \psi(y) \right) \varphi(y)
dy\\
&=  \lim_{t \downarrow 0} \int_\mathcal{O} \varphi(y) \left(
\frac1t \int_\mathcal{O} \psi(x) ( p_t(x,y) dx - \delta_y(dx) )
\right) dy\\
&= \int_\mathcal{O} \varphi(y) \left( \int_\mathcal{O} \psi(x)
\frac{\partial}{\partial t}p_t(x,y) \Big|_{t=0} dx \right) dy
\end{aligned}
$$
and since the transition density function satisfies backward
parabolic equation, it follows that
$$
\langle \mathcal{A} \varphi, \psi \rangle = \int_\mathcal{O}
\varphi(y) \left( \int_\mathcal{O} \psi(x) A_x p_t(x,y)
\Big|_{t=0} dx \right) dy,
$$
where subscript $x$ denotes that the operator $A$ acts on
$p_t(x,y)$ as a function of $x$ with $t$ and $y$ fixed. Thus we
have
$$
\begin{aligned}
\langle \mathcal{A} \varphi, \psi \rangle  &= \int_\mathcal{O}
\varphi(y) \left( \int_\mathcal{O} \psi(x) A_x \delta_y(dx)
\right) dy\\
&= \int_\mathcal{O} \varphi(y) \langle A_x \delta_y, \psi \rangle
dy = \int_\mathcal{O} \varphi(y) \langle \delta_y, A_x^* \psi
\rangle
dy\\
&= \int_\mathcal{O} \varphi(y) A_y^* \psi(y) dy = \langle \varphi,
A^* \psi \rangle.
\end{aligned}
$$
Thus $\varphi$ is a weak solution to $\mathcal{A} \varphi = A
\varphi$. By Lemma \ref{l-AH}, $\varphi \in C^2(\mathcal{O})$ and
$\varphi$ is a strong solution to $\mathcal{A} \varphi = A
\varphi$. Hence $\mathcal{A} \varphi = A \varphi$ and $A \varphi
\in E$.

\smallskip \textbf{Step 2.} We will show that $\mathcal{E} \subset
\mathcal{D}(\mathcal{A})$. Let $\varphi \in \mathcal{E}$. Then
from It\^o's formula
$$
\varphi(r_{t \wedge T_n}) \e^{\alpha \int_0^{t \wedge T_n} r_s ds}
= \varphi(r) + \int_0^{t \wedge T_n} \e^{\alpha \int_0^s r_u du} A
\varphi(r_s) ds + M_{t \wedge T_n},
$$
where $T_n = \inf \{ t \geq 0 \colon |r_t| \geq n \}$ and
$$
M_{t \wedge T_n} = \int_0^{t \wedge T_n} \sigma(r_s) \varphi'(r_s)
\e^{\alpha \int_0^s r_u du} dW_s
$$
is a martingale. Taking expectations and next passing to limit
with $n$, thanks condition \textbf{(H.1)}, we obtain
$$
\mathbb{E}^r \varphi(r_t) \e^{\alpha \int_0^t r_s ds} = \varphi(r)
+ \int_0^t \mathbb{E}^r \e^{\alpha \int_0^s r_u du} A \varphi(r_s)
ds,
$$
which means that $P_t \varphi(r) = \varphi(r) + \int_0^t P_s A
\varphi(r) ds$. Therefore by the mean-value theorem
$$
\lim_{t \downarrow 0} \frac{P_t \varphi(r) - \varphi(r)}{t} =
\lim_{t \downarrow 0} \frac1t \int_0^t P_s A \varphi(r) ds = A
\varphi(r),
$$
which means that $\varphi \in \mathcal{D}(\mathcal{A})$ and
$\mathcal{A} \varphi = A \varphi$. $_\square$

\section{Lipschitz modification of the HJB equation}

In this section we will find a twice continuously differentiable solution to
the equation
\begin{equation}\label{hjb-bound}
Q K(r) + (\alpha r - \gamma) K(r) + F_m(K(r)) = 0,\qquad r\in \mathcal{O},
\end{equation}
where $F_m$ is given by \eqref{Fm}.

\begin{remark}\label{rem-cm}
It is easy to verify that $F_m$ is a continuous Lipschitz function
with Lipschitz constant $L_m = \alpha m$. Moreover, $F_m \in
C^1((0,\infty))$. Equation \eqref{hjb-bound} may be interpreted as
HJB equation for Problems A and B with assumption that $C_t = c_t
V_t$ and $c_t \in [0, m]$. Hence, a solution to \eqref{K-hjb}
should be a limit of the sequence of solutions to
\eqref{hjb-bound} as $m \to \infty$.
\end{remark}

Define  $A_\gamma: = (A-\gamma)$ and  $\mathcal{A}_\gamma:=
(\mathcal{A}-\gamma)$. Note that \eqref{hjb-bound} can be written
as
\begin{equation}\label{sol-m}
-A_\gamma K = F_m(K) \qquad \textrm{in} \; \mathcal{O}.
\end{equation}
\begin{df}\label{Def}
We call $u \in C^2(\mathcal {O})$ a \textit{subsolution} to
\eqref{sol-m} if
$$
-A_\gamma u \leq F_m(u) \qquad \textrm{in} \; \mathcal{O}.
$$
We call $u$ a \textit{supersolution} if
$$
-A_\gamma u \geq F_m(u) \qquad \textrm{in} \; \mathcal{O}.
$$
\end{df}

\begin{remark}\label{RSuper}
It is easy to verify that $\underline{K} \equiv 0$ is a
subsolution to \eqref{sol-m}. Note that $\overline{K}(r) =
N(r)^{1-\alpha}$ is a supersolution, since, by Proposition
\ref{prop-NinC2}, $N^{1-\alpha} \in C^2(\mathcal{O})$ and since
$$
-A_\gamma \overline{K}(r) = F(\overline{K}(r)) + \frac{\alpha
\sigma^2(r) (\overline{K}'(r))^2}{2(1-\alpha)\overline{K}(r)} \geq
F(\overline{K}(r)) \geq F_m(\overline{K}(r)).
$$
Furthermore, by \textbf{(H.4)}, $N^{1-\alpha} \in
\mathcal{D}(\mathcal{A})$ and $\Phi_N(r,v) = N^{1-\alpha}(r)
v^\alpha$ satisfies \eqref{Condition}.
\end{remark}

\begin{theorem}\label{t-sol-lip}
Let $\underline{K}^m \in \mathcal{D}(\mathcal{A})$ and
$\overline{K}^m \in \mathcal{D}(\mathcal{A})$ be a subsolution and
a supersolution to \eqref{sol-m}, respectively. Assume that
$\underline{K}^m \leq \overline{K}^m$. Define $K_0^m =
\underline{K}^m$ and $K_{n+1}^m$ as
\begin{equation}\label{n+1}
K_{n+1}^m = (\lambda_m - \mathcal{A}_\gamma)^{-1} (F_m(K_n^m) +
\lambda_m K_n^m),
\end{equation}
where $\lambda_m$ is such that \eqref{F+l} holds. Then $K^m$
defined as a pointwise limit of $\{K_n^m\}$, i.e.
\begin{equation}\label{u-m}
K^m(r) = \lim_{n \to \infty} K_n^m(r), \qquad \forall r \in \mathcal{O},
\end{equation}
belongs to $C^2(\mathcal {O})$ and is a strong  solution to \eqref{sol-m}.
Moreover,   $\underline{K}^m \leq K^m
\leq \overline{K}^m$ for all $m$.
\end{theorem}

\noindent {\bf Proof} From
$$
-\mathcal{A}_\gamma \underline{K}^m + \lambda_m \underline{K}^m
\leq F_m(\underline{K}^m) + \lambda_m \underline{K}^m =
-\mathcal{A}_\gamma K_1^m + \lambda_m K_1^m
$$
we get $(\lambda_m - \mathcal{A}_\gamma)(K_1^m - \underline{K}^m)
\geq 0$. Since $P_t \varphi \geq 0$ and consequently $(\lambda_m -
\mathcal{A}_\gamma)^{-1} \varphi \geq 0$ for every $\varphi \geq
0$. It follows that $\underline{K}^m \leq K_1^m$.

Now we show that $K_n^m$ is a subsolution. From \eqref{F+l} and
\eqref{n+1} we have
$$
F_m(K_1^m) + \lambda_m K_1^m \geq F_m(\underline{K}^m) + \lambda_m
\underline{K}^m = -\mathcal{A}_\gamma K_1^m + \lambda_m K_1^m,
$$
which, with help of Lemma \ref{l-Generator}, implies that $K_1^m$
is a subsolution to \eqref{sol-m}. Hence, by induction, $K_n^m
\leq K_{n+1}^m$ and $K_{n+1}^m$ is a subsolution for all $n \in
\mathbb{N}_0$.

Now we show, by induction, that $K_n^m \leq \overline{K}^m$ for
all $n$. By definition $K_0^m \leq \overline{K}^m$. Assume that
$K_n^m \leq \overline{K}^m$. Then from \eqref{n+1} and \eqref{F+l}
we have
$$
-\mathcal{A}_\gamma K_{n+1}^m + \lambda_m K_{n+1}^m = F_m(K_n^m) +
\lambda_m K_n^m \leq F_m(\overline{K}^m) + \lambda_m
\overline{K}^m.
$$
Hence,
$$
-\mathcal{A}_\gamma K_{n+1}^m + \lambda_m K_{n+1}^m \leq
F_m(\overline{K}^m) + \lambda_m \overline{K}^m \leq
-\mathcal{A}_\gamma \overline{K}^m + \lambda_m \overline{K}^m
$$
implies that $(\lambda_m - \mathcal{A}_\gamma)(\overline{K}^m -
K_{n+1}^m) \geq 0$, and we obtain $K_{n+1}^m \leq \overline{K}^m$.

Summing up, we have
$$
\underline{K}^m \leq K_1^m \leq K_2^m \leq \ldots \leq K_n^m \leq
\ldots \leq \overline{K}^m \qquad \textrm{in} \; \mathcal{O}.
$$
Therefore $K^m(r)$ given by \eqref{u-m} exists for all $r$. Since
$F_m$ is continuous,
$$
F_m(K^m(r)) = \lim_{n \to \infty} F_m(K_n^m(r)),  \qquad \forall \, r
\in \mathcal{O},
$$
and from \eqref{n+1} we have
$$
\int_{\mathcal{O}} \varphi(r)(\lambda_m -\mathcal{A}_\gamma) K_{n+1}^m(r) dr
= \int_{\mathcal{O}} (F_m(K_n^m(r)) + \lambda_m K_n^m(r)) \varphi(r) dr,
$$
for any test function  $\varphi \in C_0^\infty(\mathcal{O})$. By
Lemma \ref{l-Generator},
$$
\int_{\mathcal{O}} K_{n+1}^m(r) (\lambda_m -A_\gamma^*) \varphi(r) dr =
\int_{\mathcal{O}} (F_m(K_n^m(r)) + \lambda_m K_n^m(r)) \varphi(r) dr.
$$
Let $n \to \infty$. By the dominated convergence theorem, we get
\begin{equation}\label{sol-weak-m}
- \int_{\mathcal{O}} K^m(r) A_\gamma^* \varphi(r) dr =
\int_{\mathcal{O}} F_m(K^m(r))
\varphi(r) dr.
\end{equation}
Since $K^m \leq \overline{K}^m$ and $\overline{K}^m$ is
continuous, then $K^m$ is locally bounded. Hence, $K^m$ is a weak
solution to \eqref{sol-m}, and we conclude by Lemma \ref{l-AH}.
$_\square$

\section{Proof of Theorem \ref{t-sol-a}}\label{sec-proofA}
Let $\{K^m\}$ be the sequence constructed in the previous section.
By \eqref{boundary} and \eqref{bnd-uni-int} the function
$\Phi_N(r,v)= N^{1-\alpha}(r) v^\alpha$ satisfies
\eqref{Condition}. So does $\Phi_0(r,v) \equiv 0$. Hence, Remark
\ref{RSuper} and Theorem \ref{t-sol-lip} guarantee that
$\Phi_m(r,v)= K^m(r) v^\alpha$ satisfies \eqref{Condition}.
Therefore, by Remark \ref{rem-cm}, $\Phi_m(r,v)$ is the value
function for Problem A with constraint $C_t \leq m V_t$. Hence,
$\{K^m\}$ is a non-decreasing sequence and the function
$$
K(r) = \lim_{m \to \infty} K^m(r), \qquad r \in \mathcal{O},
$$
is well defined. Note that $K>0$ in $\mathcal{O}$. Indeed, from
the continuity of $r_t$ we have $\e^{\int_0^t r_s ds}>0$,
$\mathbb{P}$-a.s. for all $t>0$ and $r \in \mathcal{O}$, which
implies  $\mathbb{E}^r \e^{\int_0^t r_s ds}>0$, and therefore
$$
K_1^1(r) = (\lambda_1 - \mathcal{A}_\gamma)^{-1} 1 = \int_0^\infty
\e^{-(\lambda_1+\gamma)t} \mathbb{E}^r \e^{\alpha \int_0^t r_s ds}dt
$$
is strictly positive in $\mathcal {O}$. Since $K_1^1 \leq K^1 \leq
K^2 \leq \ldots \leq K^m \leq \ldots \leq K$, we have $K>0$. Thus,
in particular, $F(K)$ is well defined, where $F(y) =
(1-\alpha)y^{\frac{\alpha}{\alpha-1}}$, for every $y>0$.

We will show that $K$ is a weak solution to
\begin{equation}\label{sol}
-A_\gamma K = F(K) \qquad \textrm{in} \; \mathcal{O}.
\end{equation}

To do this define
$$
Z_m = \{r \in \mathcal{O}\colon K^m(r) \geq m^{\alpha-1} \}.
$$
Clearly $Z_m \subset Z_{m+1}$ for all $m \in \mathbb{N}$. Since
$F_m(y) = F(y)$ for every $y \geq m^{\alpha -1}$, we have
$$
F_m(K^m(r)) = F(K^m(r)), \qquad \forall\,  r \in Z_n \ \forall \, m
\geq n,
$$
which implies, from continuity of $F$, that for any $r \in
\bigcup_{n=1}^\infty Z_n$,
\begin{equation}\label{F} \lim_{m \to
\infty} F_m(K^m(r)) = \lim_{m \to \infty} F(K^m(r)) = F(K(r)).
\end{equation}

Now we show that \eqref{F} holds for any $r \in \mathcal{O}$. To
do this note that $\bigcup_{n=1}^\infty Z_n = \mathcal{O}$.
Indeed, since $K^1(r)
>0$ for any $r \in \mathcal{O}$, then for any $r \in \mathcal{O}$ there is such $m$, that
$$
\quad K^m(r) \geq K^1(r) > m^{\alpha-1} > 0,
$$
and hence  $r \in Z_m$.

Note that for any $m>1$ and $r \in \mathcal{O}$ we have
$$
|F_m(K^m(r)) \varphi(r)| \leq (1-\alpha)
(K^1(r))^\frac{\alpha}{\alpha-1} |\varphi(r)|.
$$
Let $m \to \infty$ in \eqref{sol-weak-m}. By the inequality above
and the dominated convergence theorem, we get
$$
- \int_{\mathcal{O}} K(r) A_\gamma^* \varphi(r) dr = \int_{\mathcal{O}} F(K(r)) \varphi(r) dr, \qquad \forall\, \varphi \in C_0^\infty (\mathcal{O}),
$$
which means that $K$ is a weak solution to \eqref{sol}, whenever
$K$ is locally integrable. Since $\overline{K}^m = N^{1-\alpha}
\in E$, we have $K \leq N^{1-\alpha}$ and from continuity of
$N^{1-\alpha}$, the  function $K$ is locally bounded. By Lemma
\ref{l-AH}, $K$ is a strong solution to \eqref{sol}.

By \eqref{boundary} and \eqref{bnd-uni-int}, $\Phi(r,v)= K(r)
v^\alpha$ satisfies the boundary condition \eqref{Condition}.

\section{Solution to Problem B}

This section provides the existence and approximating scheme for
the solution $K$ to the HJB equation \eqref{K-hjb} for Problem B.
Let $(\tilde{E}, \|\cdot\|_{\tilde{E}})$ be a Banach space of
continuous functions on $\mathcal{O}^+$.

Recall that $\tau_0^r =\inf \{ t \geq 0 \colon r_t = 0 \}$. For
all $r \in \mathcal{O}^+$ we define the following functions:
$$
\tilde{N}(r) = \mathbb{E}^r
\int_0^{\tau_0^r}\e^{\frac{1}{1-\alpha}(-\gamma t + \alpha
\int_0^t r_s ds)}dt
$$
and
$$
K_L(r) = \mathbb{E}^r \e^{-\gamma \tau_0^r + \alpha
\int_0^{\tau_0^r} r_s ds}.
$$

Let $K_U(r) = K_L(r) + \tilde{N}^{1-\alpha}(r)$ for all $r \in
\mathcal{O}^+$ and let $(\tilde{r}_t) = (r_{t \wedge \tau_0^r})$.
We denote by $\mathbf{(\widetilde{H.1})}$,
$\mathbf{(\widetilde{H.2})}$ and $\mathbf{(\widetilde{H.3})}$ the
equivalents to \textbf{(H.1)}, \textbf{(H.2)} and \textbf{(H.3)}
respectively, where $r \in \mathcal{O}^+$, and $(r_t)$ and $E$ are
replaced by $(\tilde r_t)$ and $\tilde E$.

Clearly, $K_L \leq K_U$, and the following hypothesis is needed to
show that the boundary condition \eqref{Condition-B} holds for any
continuous function $f$ satisfying $K_L \leq f \leq K_U$ in
$\mathcal{O}^+$.

\smallskip \noindent $\mathbf{(\widetilde{H.4})}$
For any $r \in \mathcal{O}^+$,
\begin{equation} \label{boundary-B}
\lim_{t \to \infty} \mathbb{E}^r \e^{-\gamma t + \alpha \int_0^t
r_s ds} K_U(r_t) \chi_{\{\tau_B^{(C;r,v)} = \infty\}} = 0
\end{equation}
whenever $\mathbb{P}^r(\tau_B^{(C;r,v)} = \infty)>0$, and for any
stopping time $\tau_B^{(C;r,v)}$,
\begin{equation} \label{bnd-uni-int-B}
\left\{ \e^{-\gamma\tau_n + \alpha \int_0^{\tau_n} r_s ds}
K_U(r_{\tau_n}) \chi_{\{\tau_B^{(C;r,v)} < \infty\}} \right\}_{n
\in \mathbb{N}}
\end{equation}
is uniformly integrable, where $\tau_n = n \wedge
\tau_B^{(C;r,v)}$. Moreover, $K_L\in \mathcal{D}(\mathcal{A})$ and
$K_U\in \mathcal{D}(\mathcal{A})$, where $\mathcal{D}(\mathcal{A})
= \{ \varphi \in C^2(\mathcal{O}^{++}) \cap \tilde{E} \colon  A
\varphi \in \tilde{E}\}$, and $K_L(0) = K_U(0) = 1$.

\smallskip
Note that if $\mathbf{(\widetilde{H.4})}$ holds, then it holds
simultaneously for both processes $(r_t)$ and $(\tilde{r}_t)$.

\smallskip
Assume additionally

\smallskip \noindent $\mathbf{(\widetilde{H.5})}$  For any $r \in \mathcal{O}^{++}$, one has $A_\gamma K_L(r) =0$.

\smallskip
By ${\bf (\widetilde{H.5})}$, $K_L$ is a subsolution to
\eqref{hjb-bound}. It is easy to see that under assumptions of
Proposition \ref{prop-NinC2}, $\tilde{N}$ satisfies \eqref{N-hjb};
it is enough to take $a_n = 0$ in the proof. Thus $\tilde{N}$ is a
supersolution to \eqref{hjb-bound}. Hence, by ${\bf
(\widetilde{H.5})}$,
$$
\begin{aligned}
A_\gamma K_U + F_m(K_U) &= A_\gamma \tilde{N}^{1-\alpha} + F_m(K_L
+ \tilde{N}^{1-\alpha})\\
&\leq A_\gamma \tilde{N}^{1-\alpha} + F_m(\tilde{N}^{1-\alpha})
\leq 0
\end{aligned}
$$
and $K_U$ is also a supersolution.

Since $K_L(0)=K_U(0)=1$, then from the fact that $K_L \leq K \leq
K_U$ (see Theorem \ref{t-sol-b} below) we have a condition
$K(0)=1$, which with help of \eqref{boundary-B} and
\eqref{bnd-uni-int-B} implies \eqref{Condition-B}. Furthermore,
the value function $\Phi(\cdot,v) \in C^2(\mathcal{O}^{++}) \cap
C(\mathcal{O})$ for any $v > 0$. The proof of the following result
is analogous to that of Theorem \ref{t-sol-a} and is left to the
reader.

\begin{theorem}\label{t-sol-b}
Assume that $\mathbf{(\widetilde{H.1})}$ --
$\mathbf{(\widetilde{H.5})}$ are fulfilled. Then there is a
solution $K$ to \eqref{K-hjb} with condition \eqref{Condition-B}.
Moreover, $K_L(r) \leq K(r) \leq K_U(r)$, $r \in \mathcal{O}^+$.
Finally, for any sequence $\{\lambda_m\} \subset
\varrho(\mathcal{A}_\gamma)$ satisfying \eqref{F+l} for any $m >
0$, one has
$$
K(r) = \lim_{m \to \infty} \lim_{n \to \infty} K^m_n(r), \qquad r
\in \mathcal{O}^+,
$$
where $\{K^m_n\}$ is a non-decreasing sequence of both $m$ and
$n$, defined as follows
$$
\begin{aligned}
& K^m_0 = K_L,\\
& K^m_{n+1} = (\lambda_m - \mathcal{A}_\gamma)^{-1}(F_m(K^m_n) +
\lambda_m K^m_n).
\end{aligned}
$$
\end{theorem}

\section{Vasicek model}\label{SVasicek}
Let us recall that in the so-called  \textit{Vasicek model}
$(r_t)$ is given by
\begin{equation}\label{Vck-dr}
dr_t= (a - b r_t) dt + \sigma dW_t,
\end{equation}
with $a, b, \sigma >0$. Let
\begin{equation}\label{E-Vck}
E = \{\varphi \in C(\mathbb{R})\colon  \lim_{|r| \to \infty}
|\varphi(r)| \e^{-\frac\alpha b |r|} = 0 \}
\end{equation}
and
$$
\|\varphi\|_E = \sup_{r \in \mathbb{R}} |\varphi(r)|
\e^{-\frac\alpha b |r|}.
$$

\begin{theorem}\label{t-Vck}
The assumptions of Theorem \ref{t-sol-a} are satisfied, whenever
\begin{equation}\label{gam-N-Vck}
\gamma > \max \{ \gamma_1, \gamma_2 \},
\end{equation}
where
$$
\gamma_1 = \frac{\alpha a}{b} + \frac{\alpha^2
\sigma^2}{(1-\alpha) b^2} \qquad {\rm and} \qquad \gamma_2 =
\frac{\alpha a}{b} + \frac{3\alpha^2 \sigma^2}{2 \sqrt{1-\alpha}
b^2} + \alpha \sigma \frac{b+1}{b}.
$$
\end{theorem}
\noindent\textbf{Proof} Notice that for any stopping time $T_n$,
$$
\begin{aligned}
| \varphi(r_{t \wedge T_n}) | \e^{ \alpha \int_0^{t \wedge T_n}
r_s ds} &\leq \| \varphi \|_E \e^{ \frac{\alpha}{b} |r_{t
\wedge T_n}| + \alpha \int_0^{t \wedge T_n} |r_s| ds}\\
&\leq \| \varphi \|_E \e^{(\frac{\alpha}{b} + \alpha t) \sup_{0
\leq s \leq t} |r_s| }
\end{aligned}
$$
and, by Fernique's theorem, the r.h.s. is integrable for any fixed
$t \geq 0$. We similarly obtain
$$
\left| \int_0^{t \wedge T_n} \varphi(r_s) \e^{ \alpha \int_0^s r_u
du} ds \right| \leq  \| \varphi \|_E t \e^{(\frac{\alpha}{b} +
\alpha t) \sup_{0 \leq s \leq t} |r_s| }.
$$
Therefore \textbf{(H.1)} is satisfied.

It is easy to check that $(E,\|\cdot\|_E)$ satisfies
\textbf{(H.2)}. Assume that $(r_t)$ is given by \eqref{Vck-dr} and
that  $(P_t)$ is given by \eqref{P}. In Appendix $A$ it is shown
that $(P_t)$ is a $C_0$-semigroup on $E$, and hence hypothesis
\textbf{(H.3)} is satisfied. Therefore we have to show
\textbf{(H.4)}. We split a verification of \textbf{(H.4)} into
several steps.

\smallskip \textbf{Step 1.}
First we show that $N(r)<\infty$ for any $r \in \mathbb{R}$. From
\eqref{Vck-dr} we obtain
\begin{equation}\label{rX}
r_t = r \e^{-bt} + \tfrac{a}{b} (1 - \e^{-bt}) + \sigma X_t,
\end{equation}
where
\begin{equation}\label{X}
X_t = \int_0^t \e^{-b(t-s)} dW_s
\end{equation}
and its distribution does not depend on $r$. In what follows we
denote by $\mathcal{L}(\xi)$ the law (distribution) of a random
variable $\xi$. Note that
\begin{equation}\label{X-dist}
\mathcal{L}(X_t) = \mathcal{N}\left(0,\tfrac{1}{2b}(1 - \e^{-2bt})
\right).
\end{equation}
Therefore, we have
$$
\begin{aligned}
N(r) &= \mathbb{E}^r \int_0^\infty \e^{-\frac{\gamma}{1-\alpha} t
+\frac{\alpha}{1-\alpha} \int_0^t r_s ds} dt\\
&\leq \mathbb{E}^r \int_0^\infty \e^{-\frac{\gamma}{1-\alpha} t
+\frac{\alpha}{1-\alpha} \int_0^\infty |r| \e^{-bs} ds
+\frac{\alpha}{1-\alpha}
\tfrac{a}{b}t +\frac{\alpha}{1-\alpha} \sigma \int_0^t X_s ds} dt\\
&= \e^{\frac{\alpha}{(1-\alpha)b} |r|} \int_0^\infty
\e^{-\frac{\gamma b - \alpha a}{(1-\alpha)b} t} \mathbb{E}
\e^{\frac{\alpha}{1-\alpha} \sigma \int_0^t X_s ds} dt
\end{aligned}
$$
and by Fubini's theorem
\begin{equation}\label{Y}
Y_t := \int_0^t X_s ds = \int_0^t \int_u^t \e^{-b(s-u)} ds dW_u =
\frac1b \int_0^t (1- \e^{-b(t-u)}) dW_u,
\end{equation}
which implies that
\begin{equation}\label{Y-dist}
\mathcal{L}(Y_t) = \mathcal{N}\left(0,\tfrac{1}{b^2}(t-
\tfrac{3}{2b} + \tfrac{2}{b} \e^{-bt} - \tfrac{1}{2b} \e^{-2bt})
\right).
\end{equation}
Thus
$$
\begin{aligned}
N(r) &\leq \e^{\frac{\alpha}{(1-\alpha)b} |r|} \int_0^\infty
\e^{-\frac{\gamma b - \alpha a}{(1-\alpha)b} t} \e^{\frac{\alpha^2
\sigma^2}{2 (1-\alpha)^2 b^2} \left(t- \frac{3}{2b} +
\frac{2}{b} \e^{-bt} - \frac{1}{2b} \e^{-2bt} \right)} dt\\
&\leq \e^{\frac{\alpha}{(1-\alpha)b} |r|} \int_0^\infty
\e^{-\frac{\gamma b - \alpha a}{(1-\alpha)b} t} \e^{\frac{\alpha^2
\sigma^2}{2 (1-\alpha)^2 b^2} t} dt\\
&= \e^{\frac{\alpha}{(1-\alpha)b} |r|} \int_0^\infty
\e^{-\frac{1}{1-\alpha} (\gamma - \frac{\alpha a}{b} -
\frac{\alpha^2 \sigma^2}{2 (1-\alpha) b^2}) t} dt < \infty,
\end{aligned}
$$
by \eqref{gam-N-Vck}. Analogously we can show that by
\eqref{gam-N-Vck}, condition \eqref{N2} holds.

\smallskip \textbf{Step 2.} Here we show \eqref{boundary}. We have
just shown that
\begin{equation}\label{N-rho}
N(r) \leq \frac{\e^{\frac{\alpha}{(1-\alpha)b} |r|}}{\rho},
\end{equation}
where
$$
\rho = \frac{1}{1-\alpha} \left( \gamma - \frac{\alpha a}{b} -
\frac{\alpha^2 \sigma^2}{2(1-\alpha) b^2} \right)
$$
is positive by \eqref{gam-N-Vck}. Then to prove \eqref{boundary}
it is enough to show that
$$
\lim_{t \to \infty} \mathbb{E}^r \e^{-\gamma t + \alpha \int_0^t
r_s ds + \frac{\alpha}{b}|r_t|} = 0.
$$
From H\"{o}lder's inequality
$$
\overline{\lim_{t \to \infty}} \mathbb{E}^r \e^{-\gamma t + \alpha
\int_0^t r_s ds + \frac{\alpha}{b}|r_t|} \leq \overline{\lim_{t
\to \infty}} \left( \mathbb{E}^r \e^{-\frac{\gamma}{1-\alpha} t +
\frac{\alpha}{1-\alpha} \int_0^t r_s ds} \right)^{1-\alpha} \left(
\mathbb{E}^r \e^{\frac{|r_t|}{b}} \right)^{\alpha}
$$
and we easily compute that
$$
\overline{\lim_{t \to \infty}} \mathbb{E}^r
\e^{-\frac{\gamma}{1-\alpha} t + \frac{\alpha}{1-\alpha} \int_0^t
r_s ds} \leq \lim_{t \to \infty} \e^{\frac{\alpha}{(1-\alpha)b}
|r| - \rho t} = 0, \qquad \forall r \in \mathbb{R}.
$$

Note that given $\xi$ with distribution $\mathcal{N}(m,s^2)$, it
is easy to show that
\begin{equation}\label{EeNorm}
\mathbb{E} \e^{\kappa |\xi|} \leq \e^{\frac{\kappa^2
s^2}{2}}(1+\e^{\kappa |m|}).
\end{equation}
Therefore, the expression in the second bracket above is dominated
by $ \e^{\frac{\sigma^2}{4b^3}} \left(1+\e^{\frac{|r|}{b}
+\frac{a}{b^2}} \right)$, which is finite for every $r \in
\mathbb{R}$. Thus \eqref{boundary} holds.

\smallskip \textbf{Step 3.} Here we show that the family in
\eqref{bnd-uni-int} is uniformly integrable. By the de la Vall\'ee
Poussin theorem (see e.g. \cite{Oksendal}, p. 241), it is enough
to show that
\begin{equation}\label{ui}
\sup_{n \in \mathbb{N}} \mathbb{E}^r \left( \e^{-\gamma\tau_n +
\alpha \int_0^{\tau_n} r_s ds} N^{1-\alpha}(r_{\tau_n})
\right)^{\frac{1}{\beta}} < \infty, \qquad \forall r \in
\mathcal{O},
\end{equation}
for some $\beta<1$. Here we take $\beta = \sqrt{1-\alpha}$.

By \eqref{N-rho} we obtain
$$
\begin{aligned}
\e^{-\frac{\gamma}{\beta} \tau_n + \frac{\alpha}{\beta}
\int_0^{\tau_n} r_s ds} N^\beta(r_{\tau_n}) &\leq \rho^{-\beta}
\e^{-\frac{\gamma}{\beta} \tau_n + \frac{\alpha}{\beta}
\int_0^{\tau_n}| r_s| ds + \frac{\alpha}{b
\beta} |r_{\tau_n}|}\\
&\leq \rho^{-\beta} \e^{\sup_{t \leq n} (-\frac{\gamma}{\beta} t +
\frac{\alpha}{\beta} \int_0^t |r_s| ds +
\frac{\alpha}{b \beta} |r_t|)} \leq I,\\
\end{aligned}
$$
where
$$
I = \rho^{-\beta} \e^{\frac{a \alpha}{b^2 \beta} +
\frac{2\alpha}{b \beta} |r|} \e^{\sup_{t \leq n} (-\frac{b \gamma
- a \alpha}{b \beta} t + \frac{\alpha \sigma}{b \beta} |X_t| +
\frac{\alpha \sigma}{\beta} \int_0^t |X_s| ds)},
$$
and $X_t$ is defined by \eqref{X}. Let $g(x) = \sqrt{1+x^2}$ and
$h(x) = x^2/\sqrt{1+x^2}$. Then
$$
I \leq \rho^{-\beta} \e^{\frac{a \alpha}{b^2 \beta} +
\frac{2\alpha}{b \beta} |r|} \e^{\sup_{t \leq n} (-\frac{b \gamma
- a \alpha - b \alpha \sigma}{b \beta} t + \frac{\alpha \sigma}{b
\beta} g(X_t) + \frac{\alpha \sigma}{\beta} \int_0^t h(X_s) ds)}.
$$
By It\^o's formula
$$
\frac{\alpha \sigma}{b \beta} g(X_t) + \frac{\alpha \sigma}{\beta}
\int_0^t h(X_s) ds = \frac{\alpha \sigma}{b \beta} + \Psi_t + R_t,
$$
where
$$
\Psi_t = \frac{\alpha \sigma}{b \beta} \int_0^t g'(X_s) dW_s -
\frac12 \left(\frac{\alpha \sigma}{b \beta}\right)^2 \int_0^t
(g'(X_s))^2 ds
$$
and
$$
R_t = \frac12 \int_0^t \left( \frac{\alpha \sigma}{b \beta}
g''(X_s) + \left(\frac{\alpha \sigma}{b \beta} g'(X_s)\right)^2
\right) ds.
$$

Note that $|g'(x)|<1$ and $|g''(x)|<2$. Therefore by the Novikov
condition $M_t = \e^{\Psi_t}$ is a martingale, and $R_t < (
\frac{\alpha \sigma}{b \beta}  + \frac12 (\frac{\alpha \sigma}{b
\beta})^2 )t$. By \eqref{gam-N-Vck} there is a $\kappa >
(\frac{\alpha \sigma}{b \beta})^2$ such that
$$
\frac{b \gamma - a \alpha - b \alpha \sigma}{b \beta} >
\frac{\alpha \sigma}{b \beta} + \frac12 \left(\frac{\alpha
\sigma}{b \beta}\right)^2 + \kappa.
$$
Then
$$
\sup_{n \in \mathbb{N}} \mathbb{E}^r \e^{-\frac{\gamma}{\beta}
\tau_n + \frac{\alpha}{\beta} \int_0^{\tau_n} r_s ds}
N^\beta(r_{\tau_n}) \leq \rho^{-\beta} \e^{\frac{a \alpha + b
\alpha \sigma}{b^2 \beta} + \frac{2\alpha}{b \beta} |r|} \sup_{n
\in \mathbb{N}} \mathbb{E} \sup_{t \leq n} M_t \e^{-\kappa t}
$$
and it is enough to show that
$$
\sup_{n \in \mathbb{N}} \mathbb{E} \sup_{t \leq n} M_t \e^{-\kappa
t} < \infty.
$$
We have
$$
\begin{aligned}
\mathbb{E} \sup_{t \leq n} M_t \e^{-\kappa t} &\leq
\sum_{j=0}^{n-1} \mathbb{E} \sup_{t \in [j,j+1]} M_t \e^{-\kappa
t} \leq \sum_{j=0}^{n-1} \e^{-\kappa j} \mathbb{E} \sup_{t \in
[j,j+1]} M_t\\
&\leq \sum_{j=0}^{n-1} \e^{-\kappa j} (1+ \mathbb{E} (\sup_{t \leq
j+1} M_t)^2) \leq \sum_{j=0}^{n-1} \e^{-\kappa j} (1+ 4 \mathbb{E}
M_{j+1}^2),
\end{aligned}
$$
where the last estimate holds due to Doob's inequality. Since $M_t
\leq \e^{(\frac{\alpha \sigma}{b \beta})^2 t} \tilde{M}_t$, where
$\tilde{M}_t$ is a martingale of the same form as $M_t$, but with
constant $2\frac{\alpha \sigma}{b \beta}$ instead of $\frac{\alpha
\sigma}{b \beta}$, then
$$
\sup_{n \in \mathbb{N}}\mathbb{E} \sup_{t \leq n} M_t \e^{-\kappa
t} \leq \sup_{n \in \mathbb{N}} \sum_{j=0}^{n-1} \e^{-\kappa j}
(1+ 4 \e^{(\frac{\alpha \sigma}{b \beta})^2(j+1)}) < \infty.
$$

\smallskip \textbf{Step 4.} Here we show that $N^{1-\alpha} \in C^2(\mathbb{R}) \cap E$. By Proposition \ref{prop-NinC2}, $N^{1-\alpha} \in
C^2(\mathbb{R})$. To show that $N^{1-\alpha} \in E$ we have to
prove that
$$
\lim_{|r| \to +\infty} N^{1-\alpha}(r) \e^{-\frac\alpha b |r|} =
0.
$$
It is easy to see that
$$
\lim_{r \to -\infty} N^{1-\alpha}(r) \e^{-\frac\alpha b |r|} =
\lim_{r \to -\infty} N(r) \e^{-\frac{\alpha}{(1-\alpha) b} |r|} =
0.
$$
The condition
$$
\lim_{r \to +\infty} N(r) \e^{-\frac{\alpha}{(1-\alpha) b} |r|} =
0
$$
amounts to
$$
\lim_{x \to +\infty} \int_0^\infty \e^{-k t - x \e^{-t}} dt =0,
\qquad k>0,
$$
which clearly holds.

\smallskip \textbf{Step 5.} Finally, we need to show $A
N^{1-\alpha} \in E$. By the definition of $A$ (see \eqref{A}) and
the previous step of the proof, we know that $A N^{1-\alpha} \in
C(\mathbb{R})$. Thus we need to verify the condition
$$
\lim_{|r| \to \infty} |A N^{1-\alpha}(r)|
\e^{-\frac{\alpha}{b}|r|} =0.
$$

By \eqref{N-hjb} we have
$$
A N^{1-\alpha}(r) = N^{1-\alpha}(r) \left( \gamma -
\frac{1-\alpha}{N(r)} - \frac{\alpha(1-\alpha)\sigma^2}{2} \left(
\frac{N'(r)}{N(r)} \right)^2 \right).
$$
Since
$$
N(r) = \int_0^\infty \e^{r \frac{\alpha}{(1-\alpha)b}
(1-\e^{-bt})} \phi(t) dt,
$$
where $\phi$ is, by Step 1, a strictly positive integrable
function, then $N$ is positive and increasing. Furthermore,
$$
N'(r) =  \int_0^\infty \frac{\partial}{\partial r} \e^{r
\frac{\alpha}{(1-\alpha)b} (1-\e^{-bt})} \phi(t) dt \leq
\frac{\alpha}{(1-\alpha)b} N(r).
$$
Hence,
$$
\lim_{|r| \to \infty} |A N^{1-\alpha}(r)|
\e^{-\frac{\alpha}{b}|r|} \leq \lim_{|r| \to \infty}
N^{1-\alpha}(r) \e^{-\frac{\alpha}{b}|r|} \left( \gamma +
\frac{1-\alpha}{N(r)} + \frac{\alpha^3\sigma^2}{2(1-\alpha)b^2}
\right)
$$
and, since $N^{1-\alpha} \in E$, the limit above is equal to zero.
$_\square$

\bigskip
Note that the condition $\gamma > \gamma_1$ assures the finiteness
of $N(r)$ for any $r \in \mathbb{R}$ and that assumption
\eqref{N2} holds, and the condition $\gamma
> \gamma_2$ is needed for uniform integrability of the family in
\eqref{bnd-uni-int}.

\bigskip
Let $\delta > \frac{\alpha(3-\alpha)}{b(1-\alpha)}$ and let
$$
\tilde E = E_\delta = \{\varphi \in C([0,\infty))\colon  \lim_{r
\to \infty} |\varphi(r)| \e^{-\delta r} = 0 \}
$$
be equipped with the norm
$$
\|\varphi\|_{\tilde E} = \sup_{r \in [0,\infty)} |\varphi(r)|
\e^{-\delta r}.
$$
Then we have the following result.

\begin{theorem}\label{t-Vck-b}
The assumptions of Theorem \ref{t-sol-b} are fulfilled in the
Vasicek model \eqref{Vck-dr}, whenever \eqref{gam-N-Vck} holds.
\end{theorem}

\noindent\textbf{Proof} Verification of
$\mathbf{(\widetilde{H.1})}$, $\mathbf{(\widetilde{H.2})}$ and
$\mathbf{(\widetilde{H.3})}$ is left to the reader, as it is
similar to verification of \textbf{(H.1)}, \textbf{(H.2)} and
\textbf{(H.3)}. Here we verify only $\mathbf{(\widetilde{H.4})}$
and $\mathbf{(\widetilde{H.5})}$. To show
$\mathbf{(\widetilde{H.5})}$ define a sequence of functions
$\{K_L^n\}$, such that
$$
\left\{ \begin{array}{l}
  A_\gamma K_L^n(r) = 0, \quad r \in (0,n), \\
  K_L^n(0) = K_L^n(n) = 1. \\
\end{array} \right.
$$
Then, by \eqref{gam-N-Vck},
$$
K_L^n(r) = \mathbb{E}^r \e^{-\gamma \tau_{0,n}^r + \alpha
\int_0^{\tau_{0,n}^r} r_s ds},
$$
where $\tau_{0,n}^r = \tau_0^r \wedge \tau_n^r$ and $\tau_n^r =
\inf\left\{t\ge 0: r_t = n\right\}$. Furthermore, for any test
function $\varphi \in C_0^\infty((0,n))$,
$$
\int_0^\infty \varphi(r) A_\gamma K_L^n(r) dr = 0,
$$
which implies that
$$
\lim_{n\to\infty} \int_0^\infty K_L^n(r) A_\gamma^* \varphi(r) dr
= 0.
$$
Under the following conditions,
$$
\begin{aligned}
(i) &\quad \lim_{n\to\infty} \tau_{0,n}^r = \tau_0^r, \quad
\forall r>0, \forall \omega \in \Omega,\\
(ii) &\quad \mathbb{P}^r(\tau_0^r<\infty)=1, \quad \forall r>0,\\
(iii) &\quad \sup_{r \leq j} \sup_{n\in\mathbb{N}} |K_L^n(r)| <
\infty, \quad \forall j>0,
\end{aligned}
$$
which we will verify below, we have $\lim_{n\to\infty} K_L^n(r) =
K_L(r)$ and the convergence is almost uniform. Thus,
$$
\int_0^\infty K_L(r) A_\gamma^* \varphi(r) dr = 0,
$$
and it holds for any $\varphi \in C_0^\infty((0,\infty))$. By
Lemma \ref{l-AH}, $A_\gamma K_L(r) =0$.

Since any Ornstein--Uhlenbeck process is recurrent, then $(ii)$
holds. Condition $(i)$ is implied by the continuity of
trajectories. To show $(iii)$, note that due to Step 3 of the
proof of Theorem \ref{t-Vck} the sequence of random variables in
the second expression below is uniformly integrable, which
justifies the first equality, and the constant $d<\infty$ does not
depend on $r$. Thus,
$$
\begin{aligned}
\sup_{r \leq j} \sup_{n\in\mathbb{N}} |K_L^n(r)| &= \sup_{r \leq
j} \sup_{n\in\mathbb{N}} \mathbb{E}^r \lim_{m \to \infty}
\e^{-\gamma (\tau_{0,n}^r \wedge m) + \alpha \int_0^{\tau_{0,n}^r
\wedge m}r_s
ds}\\
&\leq \sup_{r \leq j} \sup_{n\in\mathbb{N}} \e^{\frac\alpha b r}
\sup_{m \geq 0} \mathbb{E} \e^{\sup_{t \leq m}(-\gamma t +
\frac{\alpha a}{b} t + \alpha \sigma Y_t)}\\
&\leq d \sup_{r \leq j} \e^{\frac\alpha b r} = d \e^{\frac\alpha b
j} < \infty,
\end{aligned}
$$
where $Y_t$ is given by \eqref{Y}.

We proceed to show that ${\bf (\widetilde{H.4})}$ holds. Since
condition $(ii)$ above holds, then we do not have to verify
\eqref{boundary-B}.

Note that $K_U \leq K_L + N^{1-\alpha}$. Therefore, by Theorem
\ref{t-Vck}, the sequence in \eqref{bnd-uni-int-B} is uniformly
integrable whenever uniformly integrable is the sequence
$$
\left\{ \e^{-\gamma\tau_n + \alpha \int_0^{\tau_n} \tilde r_s ds}
K_L(\tilde r_{\tau_n}) \right\}_{n \in \mathbb{N}}.
$$
By the strong Markov property it is equivalent to the uniform
integrability of
$$
\left\{ \mathbb{E}^r \Big[ \e^{-\gamma\tau_0^r + \alpha
\int_0^{\tau_0^r} \tilde r_s ds} \big| \mathcal{F}_n \Big]
\right\}_{n \in \mathbb{N}},
$$
which is fulfilled whenever
$$
\mathbb{E}^r \e^{-\gamma\tau_0^r + \alpha \int_0^{\tau_0^r} \tilde
r_s ds} <\infty.
$$
This holds, since $\tau_0^r \wedge n \to \tau_0^r$ as $n \to
\infty$ and
$$
\left\{ \e^{-\gamma(\tau_0^r \wedge n) + \alpha \int_0^{\tau_0^r
\wedge n} |r_s| ds} \right\}_{n \in \mathbb{N}}
$$
is uniformly integrable (see Step 3 of the proof of Theorem
\ref{t-Vck}).

Here we will show that $K_L \in \mathcal{D}(\mathcal{A})$. Since
$A_\gamma K_L =0$ then $K_L \in C^2((0,\infty))$ and it is enough
to show that $K_L \in \tilde E$. By the similar argumentation to
that in verification of $(iii)$ in the previous step of the proof,
we have
$$
\begin{aligned}
\lim_{r \to \infty} |K_L(r)| \e^{-\delta r} &= \lim_{r \to \infty}
\mathbb{E}^r \lim_{m \to \infty} \e^{-\gamma (\tau_0^r \wedge m) +
\alpha \int_0^{\tau_0^r \wedge m}r_s ds} \e^{- \delta r}\\
&\leq \lim_{r \to \infty} \e^{\frac\alpha b r - \delta r} \sup_{m
\geq 0} \mathbb{E} \e^{\sup_{t \leq m}(-\gamma t + \frac{\alpha
a}{b} t + \alpha \sigma Y_t)}\\
&\leq d \lim_{r \to \infty} \e^{\frac\alpha b r - \delta r} = 0.
\end{aligned}
$$

Now we will show that $\tilde N \in \mathcal{D}(\mathcal{A})$,
which will imply that $K_U \in \mathcal{D}(\mathcal{A})$. Since,
by condition \eqref{gam-N-Vck}, $\tilde{N}$ satisfies
\begin{equation} \label{Nt-hjb}
Q \tilde N(r) + \frac{\alpha r -\gamma}{1-\alpha} \tilde N(r) =
-1,
\end{equation}
then $\tilde N \in C^2((0,\infty))$, and consequently $\tilde
N^{1-\alpha} \in C^2((0,\infty))$. Furthermore, $\tilde
N^{1-\alpha} \in \tilde E$, since
$$
\lim_{r \to \infty} \tilde N^{1-\alpha}(r) \e^{-\delta r} \leq
\lim_{r \to \infty} N^{1-\alpha}(r) \e^{-\delta r}\\
\leq \rho^{\alpha-1} \lim_{r \to \infty} \e^{\frac\alpha b r -
\delta r} = 0.
$$
Analogously, $\tilde N^{1-\alpha} \in E_{\delta_1}$ for any
$\delta_1 > \alpha/b$, and consequently $\tilde N \in
E_{\delta_2}$ for any $\delta_2 > \frac{\alpha}{b(1-\alpha)}$.

In order to prove that $A \tilde N^{1-\alpha} \in \tilde E$, note
that by \eqref{Nt-hjb}, we have
$$
A \tilde N^{1-\alpha}(r) = \tilde N^{1-\alpha}(r) \left( \gamma -
\frac{1-\alpha}{\tilde N(r)} - \frac{\alpha(1-\alpha)\sigma^2}{2}
\left( \frac{\tilde N'(r)}{\tilde N(r)} \right)^2 \right).
$$
We need the following result.

\begin{lemma}
Let $f \in C([0,\infty)) \cap C^1((0,\infty))$ and $f' \in
E_{\delta^1}$. Then $f \in E_{\delta^2}$ for any $\delta^2 >
\delta^1$.
\end{lemma}

\noindent\textbf{Proof} Since for any $x>0$, one has $f(x) = f(0)
+ \int_0^x f'(y) dy$, then
$$
\begin{aligned}
|f(x)| \e^{-\delta^2 x} &\leq |f(0)| \e^{-\delta^2 x} +
\e^{-(\delta^2-\delta^1) x} \int_0^x |f'(y)| \e^{-\delta^1 y}
\e^{-\delta^1 (x-y)}dy\\
&\leq |f(0)| \e^{-\delta^2 x} + \frac{M}{\delta^1}
\e^{-(\delta^2-\delta^1) x},
\end{aligned}
$$
where $M = \sup_{x \geq 0} |f'(x)| \e^{-\delta^1 x}$. Since $f'
\in E_{\delta^1}$, then $M$ is finite. $_\square$

\bigskip
Going back to the proof of Theorem \ref{t-Vck-b} note that we can
rewrite \eqref{Nt-hjb} as follows
$$
\left(\frac12 \sigma^2 \tilde N' + (a-br) \tilde N \right)' = -1 -
\left(b + \frac{\alpha r -\gamma}{1-\alpha}\right) \tilde N,
$$
which implies that the l.h.s. belongs to $E_{\delta_2}$, and by
the lemma above, $\frac12 \sigma^2 \tilde N' + (a-br) \tilde N \in
E_{\delta_3}$ for any $\delta_3 > \delta_2$. Since $(a-br) \tilde
N \in E_{\delta_2} \subset E_{\delta_3}$, then $\tilde N' \in
E_{\delta_3}$. Hence, for any $\delta_1 > \alpha/b$ and $\delta =
\delta_1 + 2 \delta_3$, we get
$$
\begin{aligned}
\lim_{r \to \infty} |A_\gamma \tilde N^{1-\alpha}(r)| \e^{-\delta
r} &\leq \lim_{r \to \infty} \tilde N^{1-\alpha}(r) \e^{-\delta_1
r} \left( \gamma + \frac{1-\alpha}{\tilde N(r)} \right)
\e^{-2 \delta_3 r}\\
&\quad + \lim_{r \to \infty} \tilde N^{1-\alpha}(r) \e^{-\delta_1
r} \frac{\alpha(1-\alpha)\sigma^2}{2 \tilde N^2(r)} (\tilde N'(r)
e^{-\delta_3 r})^2.
\end{aligned}
$$
It is easy to verify that $\tilde{N}$ is increasing, which implies
that $\tilde N(r)>0$ for any $r>0$. Thus, the limit above is equal
to zero. $_\square$

\section{Invariant interval model}\label{exa-bnd}
Here we assume that the short-rate dynamics is given by \eqref{r},
$\mathcal{O}=(a, b)$, where $-\infty < a < b < \gamma/ \alpha$ and
$E = UC((a,b))$ is equipped with the supremum norm.

The sufficient condition for interval invariance is (see
\cite{Ikeda-Watanabe})
$$
s(a^+) = -\infty \qquad {\rm and} \qquad s(b^-) = \infty,
$$
where
$$
s(x) = \int_w^x \e^{\int_w^y \frac{2 \mu(z)}{\sigma^2(z)} dz}dy
$$
for a fixed $w \in (a,b)$.

It is easy to show that the conditions above holds in the model
\begin{equation}\label{dr-bnd}
dr_t = \kappa (\frac{a+b}{2}-r_t)dt + \sigma (r_t-a)(b-r_t)dW_t
\end{equation}
with $\kappa, \sigma > 0$.

\begin{theorem}\label{t-bnd}
The assumptions of Theorem \ref{t-sol-a} hold in the invariant
interval model \eqref{dr-bnd}.
\end{theorem}

\noindent \textbf{Proof} Notice that for any stopping time $T_n$
and any fixed $t \geq 0$,
$$
|\varphi(r_{t \wedge T_n})| \e^{ \alpha \int_0^{t \wedge T_n} r_s
ds} \leq \| \varphi \|_E  \e^{ \alpha \int_0^{t \wedge T_n} |r_s|
ds} \leq \| \varphi \|_E \e^{ \alpha t (|a| \vee |b|) } < \infty,
$$
and similarly
$$
\left| \int_0^{t \wedge T_n} \varphi(r_s) \e^{ \alpha \int_0^s r_u
du} ds \right| \leq \| \varphi \|_E t \e^{ \alpha t (|a| \vee |b|)
} < \infty,
$$
which means that \textbf{(H.1)} is satisfied. One may easily check
that the space $(E,\|\cdot\|_E)$ satisfies \textbf{(H.2)}. Assume
that $(P_t)$ is given by \eqref{P}. In Appendix $B$ it is shown
that $(P_t)$ is a $C_0$-semigroup on $E$, and hence hypothesis
\textbf{(H.3)} is satisfied. Thus we have to verify
\textbf{(H.4)}.

Given $b < \gamma / \alpha$, we have
$$
N(r) \leq \int_0^\infty \e^{\frac{1}{1-\alpha}(-\gamma + \alpha
b)t} = \frac{1-\alpha}{\gamma - \alpha b}
$$
and
$$
\lim_{t \to \infty} \mathbb{E}^r \e^{-\gamma t + \alpha \int_0^t
r_s ds} N^{1-\alpha}(r_t) \leq \lim_{t \to \infty} \e^{-(\gamma -
\alpha b)t} \left( \frac{1-\alpha}{\gamma - \alpha b}
\right)^{1-\alpha} = 0.
$$
Thus $N<\infty$ and \eqref{boundary} holds. In the same manner we
can see that \eqref{N2} holds.

Recall that \eqref{ui} implies uniform integrability of the family
in \eqref{bnd-uni-int}. Let $\beta = 1-\alpha$. Then \eqref{ui}
holds, since we have
$$
\sup_{n \in \mathbb{N}} \mathbb{E}^r \e^{-\frac{\gamma}{1-\alpha}
\tau_n + \frac{\alpha}{1-\alpha} \int_0^{\tau_n} r_s ds}
N(r_{\tau_n}) \leq \sup_{n \in \mathbb{N}} \frac{1-\alpha}{\gamma
- \alpha b} \e^{-\frac{\gamma - \alpha b}{1-\alpha}\tau_n} =
\frac{1-\alpha}{\gamma - \alpha b}.
$$

By Proposition \ref{prop-NinC2}, one has $N^{1-\alpha} \in
C^2((a,b))$. Note that $N$ is bounded, i.e.
$$
0 < \frac{1-\alpha}{\gamma - \alpha a} \leq N(r) \leq
\frac{1-\alpha}{\gamma - \alpha b} < \infty.
$$
Since $N$ is also increasing and continuous, then there exist
finite limits $N(a^+)$ and  $N(b^-)$. Thus $N^{1-\alpha} \in
UC((a,b))$. By \eqref{N-hjb}, we have
$$
A N^{1-\alpha}(r) = N^{1-\alpha}(r) \left( \gamma -
\frac{1-\alpha}{N(r)} - \frac{\alpha(1-\alpha) (\sigma (r-a) (b-r)
N'(r))^2}{2 N^2(r)} \right).
$$
Hence, $A N^{1-\alpha} \in C((a,b))$. Since $N'(r) \geq 0$ for all
$r \in (a,b)$, and
$$
N(\frac{a+b}{2}) - N(a^+) = \int_a^{\frac{a+b}{2}} N'(r) dr  =
\int_a^{\frac{a+b}{2}} (r-a) N'(r) \frac{1}{r-a} dr
$$
is finite, then necessarily $\lim_{r \to a^+} (r-a)N'(r) = 0$.
Analogously we get $\lim_{r \to b^-} (b-r)N'(r) = 0$. Thus there
exist limits $\lim_{r \to a^+} A N^{1-\alpha}(r)$ and $\lim_{r \to
b^-} A N^{1-\alpha}(r)$, which implies that $A N^{1-\alpha} \in
UC((a,b))$. $_\square$

\section{Models with infinite value function}\label{SInfinite}
We will show here that if $(r_t)$ is either a Brownian motion or a
geometric Brownian motion, then  the value function in Problem A
is infinite.

Let us observe first  that we may assume that optimal consumption is of the proportional form $C_t=c_t V_t$, where
$$
c_t = \left\{
\begin{array}{ll}
    C_t/V_t, & \hbox{$t < \tau_A$}, \\
    0, & \hbox{$t \geq \tau_A$}, \\
\end{array}
\right.
$$
for $\tau_A$ given by \eqref{tau} and $c_t$ is well defined. Also
in this case the HJB equation and the optimal consumption $\hat C$
have the form \eqref{K-hjb} and \eqref{C-opt} respectively.
Moreover,
$$
dV_t = (r_t - c_t) V_t dt,
$$
and consequently
$$
V_t = v \e^{\int_0^t (r_s - c_s)ds} > 0, \qquad \forall \, v>0,\ \forall\, t < \tau_A,
$$
which implies that
$$
\tau_A = \inf\left\{t\geq 0\colon \int_0^t c_s ds = \infty\right\}.
$$
Thus from now on, we assume that our consumption is of the
proportional form and
\begin{equation} \label{J-prop}
J_A(c;r,v):= v^\alpha \mathbb{E}^r \int_0^\infty \e ^{-\gamma t}
c^\alpha_t \e^{\alpha \int_0^t (r_s - c_s)ds} dt
\end{equation}
with $c_t=0$ for every $t\geq \tau_A$.

\begin{lemma}\label{L-r-const}
Assume $r_t=r=const$. Then:

\noindent $i)$ If $\gamma -\alpha r\le 0$, then there is a consumption
rate $C$ such that $J_A(C;r,v)=\infty$ for all $v>0$.

\noindent
$ii)$ If $\gamma -\alpha r > 0$, then
\begin{equation}\label{const}
\begin{aligned}
\Phi_A(r,v) =\left( \frac{\gamma -\alpha r}{1-\alpha}
\right)^{\alpha -1}v^\alpha ,&\qquad \hat{C}_t
=\frac{\gamma -\alpha r}{1-\alpha} V_t,\\
V_t = \e^{(r -\frac{\gamma -\alpha r}{1-\alpha})t}v
&=\e^{\frac{r-\gamma}{1-\alpha}t}v.
\end{aligned}
\end{equation}
\end{lemma}

\noindent {\bf Proof of i)}  Whenever $\gamma -\alpha r \leq 0$,
then \eqref{J-prop} gives us the claim with $c_t \leq \alpha r -
\gamma$.  $_\square$

\noindent {\bf Proof of ii)} Since now $\mu(r) = \sigma(r) = 0$, we have
$\hat{C}_t = K^{1/(\alpha-1)} V_t$ and
$$
(\alpha r-\gamma) K + (1-\alpha) K^{\frac{\alpha}{\alpha-1}}=0
$$
instead of \eqref{K-hjb}. If $\gamma -\alpha r > 0$, then
$$
K = \left( \frac{\gamma -\alpha r}{1-\alpha} \right)^{\alpha -1}.
\qquad _\square
$$

\begin{remark}
Recall that if $r_t >0$ for every $t \geq 0$, then Problem B
amounts to Problem A. Since condition $\gamma -\alpha r\le 0$
implies $r>0$, thus the first claim in Lemma \ref{L-r-const} holds
also for Problem B. The second claim is true for Problem B,
whenever $r>0$. Otherwise $\Phi_B(r,v) = v^\alpha$ with
$\tau_B=0$.
\end{remark}

Now we formulate necessary and sufficient conditions for
finiteness of value function $\Phi_A$. Set  $\mathcal{S} = \mathcal{O}
\times (0,\infty)$.

\begin{lemma}\label{Pro-ns}
$i)$ If $\Phi_A(r,v)$ is finite for all $(r,v) \in \mathcal{S}$, then
\begin{equation}\label{nc}
\forall r \in \mathcal{O} \; \forall c > 0 \quad \mathbb{E}^r
\int_0^\infty \e^{-\gamma t +\alpha \int_0^t (r_s-c)ds} dt <
\infty.
\end{equation}
\noindent $ii)$ If the performance functional is given by
\eqref{J-prop} and
\begin{equation}\label{sc}
\exists \delta >0 \; \exists p \in (1, \tfrac{1}{\alpha}) \;
\forall r \in \mathcal{O} \quad \mathbb{E}^r \int_0^\infty \e^{-(\gamma -
\delta)q t + \alpha q \int_0^t r_s ds} dt < \infty
\end{equation}
holds with $q = p/(p-1)$, then $\Phi_A(r,v)$ is finite for all
$(r,v) \in \mathcal{S}$.
\end{lemma}

\noindent {\bf Proof of i)} Taking $c_t=c$ constant gives us the
claim. $_\square$

\noindent {\bf Proof of ii)} Set $\delta > 0$. From \eqref{J-prop} we have
$$
J_A(c;r,v) = v^\alpha  \mathbb{E}^r \int_0^\infty \e^{-(\gamma -
\delta) t + \alpha \int_0^t r_s ds} \e^{-\delta t} c_t^\alpha
\e^{-\alpha \int_0^t c_s ds} dt,
$$
and from H\"{o}lder's inequality $J_A(c;r,v)$ is dominated by
$$
v^\alpha  \mathbb{E}^r \left( \int_0^\infty \e^{-(\gamma - \delta)q
t + \alpha q \int_0^t r_s ds} dt \right)^{\frac1q} \left(
\int_0^\infty \e^{-\delta p t} c_t^{\alpha p} \e^{-\alpha p \int_0^t
c_s ds} dt \right)^{\frac1p}.
$$
From Lemma \ref{L-r-const} with $r=0$, the expression in the
second bracket above is finite for every $\delta>0$ and $p>1$ such
that $\alpha p <1$. Thus \eqref{sc} gives us the claim. $_\square$

\begin{prop}
If $(r_t)$ is a drifted Brownian motion
$$
r_t= r + \mu t + \sigma W_t,
$$
or $(r_t)$ is a geometric Brownian motion
$$
r_t= r \e^{(\mu - \frac{1}{2}\sigma^2)t + \sigma W_t},
$$
then \eqref{nc} does not hold for any $c > 0$ and consequently the
value function $\Phi_A$ for Problem $A$ is infinite.
\end{prop}
\noindent \textbf{Proof} Notice that $\mathcal{L}(\int_0^t W_s ds)
= \mathcal{N}(0, t^3/3)$, which implies
$$
\mathbb{E} (\e^{a \int_0^t W_s ds}) = \e^{\frac{a^2 t^3}{6}}.
$$
Thus if $(r_t)$ is a
drifted Brownian motion, then
$$
\begin{aligned}
\mathbb{E}^r \int_0^\infty \e^{-\gamma t +\alpha \int_0^t
(r_s-c)ds} dt &= \int_0^\infty \e^{(\alpha r -\gamma - \alpha c)t +
\frac12 \alpha \mu t^2} \mathbb{E}\, \e^{\alpha \sigma \int_0^t W_s
ds} dt\\
&= \int_0^\infty \e^{(\alpha r -\gamma - \alpha c)t + \frac12
\alpha \mu t^2 + \frac16 \alpha^2 \sigma^2 t^3} dt = \infty
\end{aligned}
$$
for all $c \geq 0$. Hence, $\Phi_A(r,v)=\infty$.

Notice that $\e^y > y$ for all $y \in \mathbb{R}$. Thus if $(r_t)$
is a geometric Brownian motion, then we have
$$
\begin{aligned}
\mathbb{E}^r \int_0^\infty \e^{-\gamma t +\alpha \int_0^t
(r_s-c)ds} dt &= \int_0^\infty \e^{-(\gamma + \alpha c)t}
\mathbb{E}\, \e^{\alpha r \int_0^t \e^{(\mu - \frac{1}{2}\sigma^2)s +
\sigma W_s}ds} dt\\
&\geq \int_0^\infty \e^{-(\gamma + \alpha c)t} \mathbb{E}\, \e^{\alpha
r \int_0^t ((\mu - \frac{1}{2}\sigma^2)s + \sigma
W_s)ds} dt\\
&= \int_0^\infty \e^{-(\gamma + \alpha c)t + \frac12 \alpha r (\mu
- \frac12\sigma^2) t^2 + \frac16 \alpha^2 r^2 \sigma^2 t^3} dt =
\infty
\end{aligned}
$$
for all $c \geq 0$. Hence, $\Phi_A(r,v)=\infty$. Moreover
$\Phi_B(r,v)=\infty$, since in this case $r_t>0$ for every $t \geq
0$. $_\square$

\section{Numerical results}

Here we present a numerical solution for a Vasicek model with
parameters $a=0.03$, $b=0.5$ and $\sigma=0.02$. We take
$\alpha=0.5$ and $\gamma = 1.5304$, which satisfies the condition
\eqref{gam-N-Vck}. Since $\gamma > \vartheta$ (see
\eqref{lambda}), we take $\lambda_m = \alpha m + 10^{-5}$.

\begin{figure}[h]
\begin{center}
  \includegraphics[width = 0.48\textwidth]{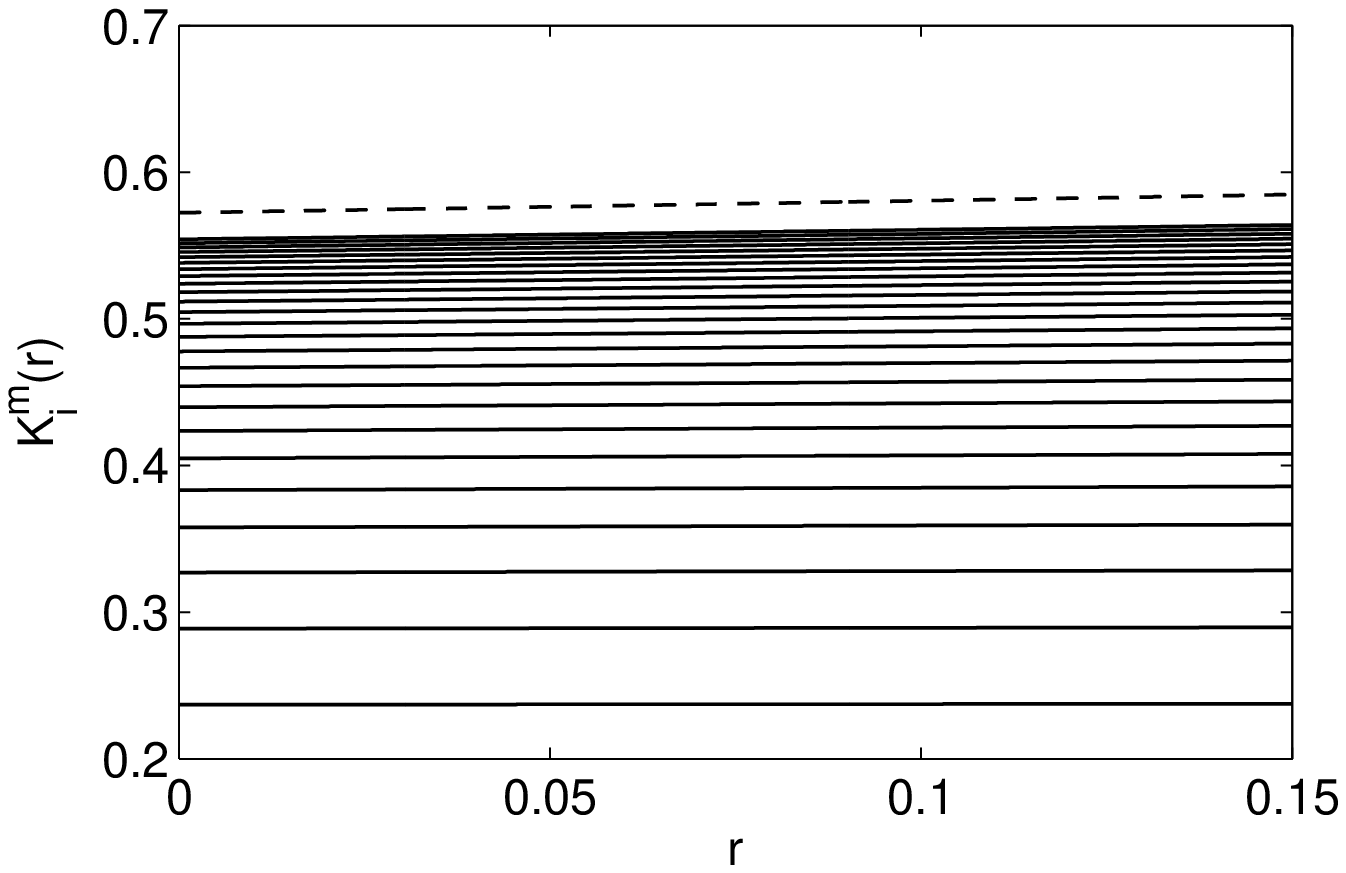} \,
  \includegraphics[width = 0.48\textwidth]{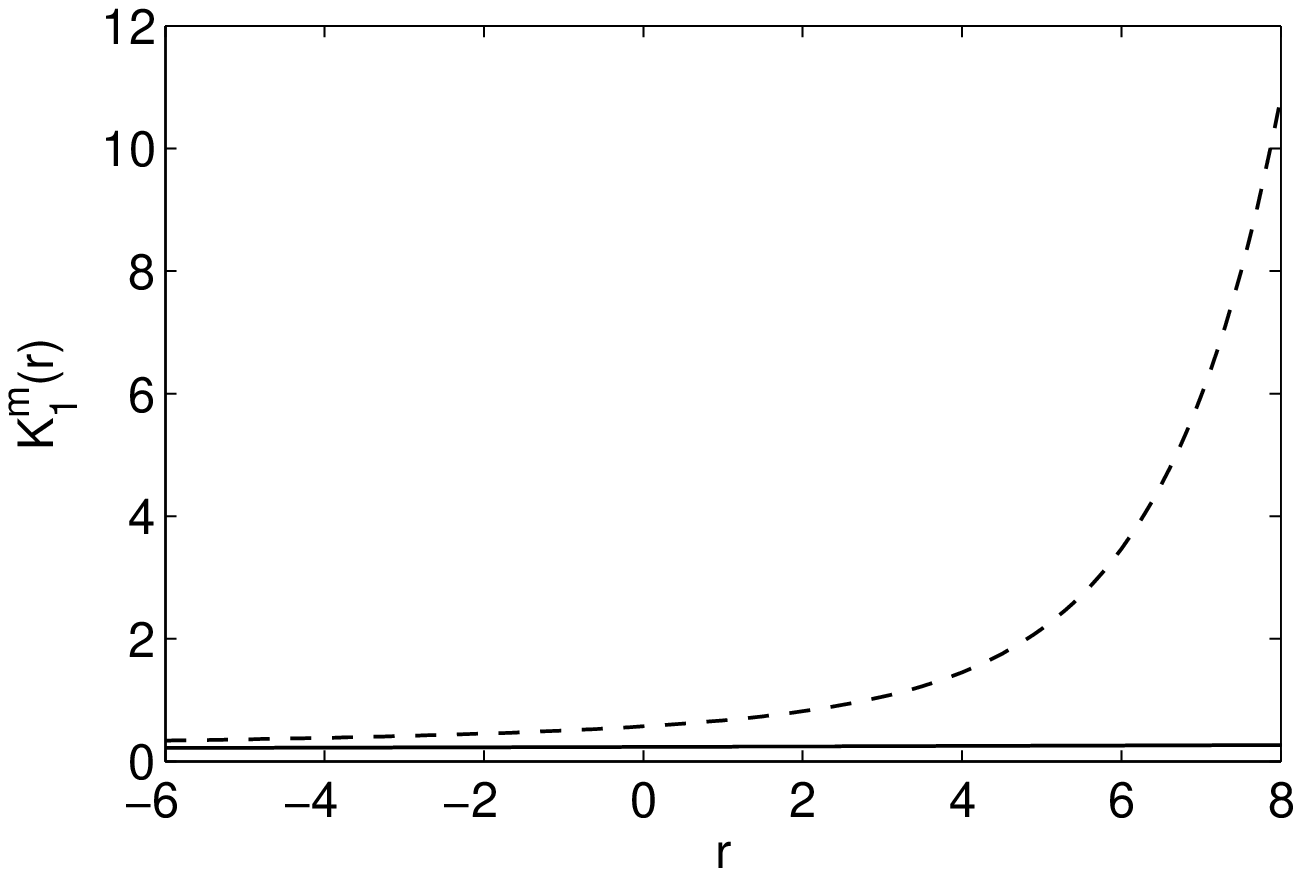}\\
  \caption{{\small (left) The dashed line is a supersolution $N^{1-\alpha}(r)$,
   and the solid lines are $K_i^m(r)$ for $i = 1,\ldots,n$. (right) The dashed line is a supersolution $N^{1-\alpha}(r)$,
   and the solid line is $K_1^m(r)$. }} \label{fig1}
\end{center}
\end{figure}

Recall that the value function is given by $\Phi(r,v) = K(r)
v^\alpha$, and $K(r)$ is as in Theorem \ref{t-sol-a}. Therefore we
have to approximate the function $K$ by $K_n^m$ for some large $m$
and $n$. Since $K_n^m(r)$ is given by recurrent formula
$$
K_n^m(r) = \int_0^\infty \int_{-\infty}^\infty S_{n-1}^m(t,r,y) dt
dy \approx \int_{t_{min}}^{t_{max}} \int_{y_{min}}^{y_{max}}
S_{n-1}^m(t,r,y) dt dy
$$
with a complicated function $S_{n-1}^m$ such that $\lim_{t
\downarrow 0} S_{n-1}^m(t,r,r) = \infty$, then we use trapezoidal
quadrature. We take $\Delta t = 0.001$ and $\Delta y = 0.0002$ to
get the result with a small error. In fact this makes the
calculations very time-consuming. Thus we take $m=65$ and $n=25$
and we have the result as in Figure \ref{fig1} for $r \in (0,
0.15)$. The result over the range $(-6,8)$ is given only for
$K_1^m$ due to very long calculations of $K_i^m$ for $i \geq 2$.

\begin{figure}[h]
\begin{center}
  \includegraphics[width = 0.9\textwidth]{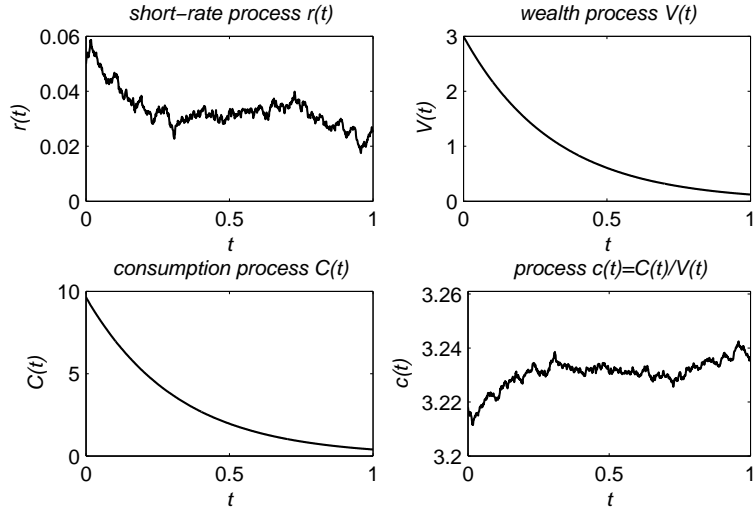}\\
  \caption{{\small Trajectories of processes $r(t,\omega)$, $V(t,\omega)$,
   $C(t,\omega)$ and $c(t,\omega)$ for the same $\omega \in \Omega$, $r
   = 0.05$ and $v=3$.}} \label{fig2}
\end{center}
\end{figure}

Next we compute trajectories of the wealth $(V_t)$, the optimal
consumption $(C_t)$ and the relative consumption $(c_t) =
(C_t/V_t)$ for a given realization of the interest rate $(r_t)$.
Clearly we take $K_n^m$ instead of $K$. The results for initial $r
= 0.05$ and $v=3$ are given in Figure \ref{fig2}.

\section*{Appendix $A$ - proof of $C_0$-semigroup property of $(P_t)$ in case of Vasicek model}

For any $\varphi \colon \mathbb{R} \mapsto \mathbb{R}$ define
$$
\|\varphi\| = \sup_{r \in \mathbb{R}} |\varphi(r)|
\e^{-\frac\alpha b |r|}.
$$
Note that $\|\varphi\| \leq \infty$, and $\|\varphi\| =
\|\varphi\|_E$ for $\varphi \in E$, where $E$ is given by
\eqref{E-Vck}. We assume that $(r_t)$ and $(P_t)$ are given by
\eqref{Vck-dr} and \eqref{P} respectively.

In the subsequent steps of the proof we need the following result.
\begin{lemma}\label{l-norm}
For any $\varphi \in E$ and any $t \geq 0$,
$$
\|P_t\varphi\| \leq 2 \e^{(\frac{\alpha^2\sigma^2}{2b^2} +
\frac{\alpha a}{b})t} \|\varphi\|.
$$
\end{lemma}
\noindent {\bf Proof} Notice that we do not assume that
$P_t\varphi \in E$. This will be shown later. Let $X_t$ and $Y_t$
be given by \eqref{X} and \eqref{Y} respectively. We have
$$
\begin{aligned}
\|P_t\varphi\| = \sup_{r \in \mathbb{R}} |P_t\varphi(r)|
\e^{-\frac\alpha b |r|} &\leq \sup_{r \in \mathbb{R}} \mathbb{E}^r
|\varphi(r_t)| \e^{\alpha\int_0^t r_s ds -\frac\alpha b |r|}\\
&\leq \|\varphi\| \sup_{r \in \mathbb{R}} \mathbb{E}^r
\e^{\alpha\int_0^t r_s ds +\frac\alpha b (|r_t|-|r|)}\\
&\leq \|\varphi\| \sup_{r \in \mathbb{R}} \e^{\frac{\alpha}{b}
(1-\e^{-bt})(r-|r|) + \frac{\alpha a}{b}t} \mathbb{E}
\e^{\alpha\sigma(\frac1 b |X_t|+Y_t)}\\
&\leq \e^{\frac{\alpha a}{b}t} \mathbb{E} \left(
\e^{\alpha\sigma(\frac1 b X_t+Y_t)} + \e^{\alpha\sigma(-\frac1 b
X_t+Y_t)} \right) \|\varphi\| \\
&\leq 2 \e^{\frac{\alpha a}{b}t}
\e^{\frac{\alpha^2\sigma^2}{2b^2}t} \|\varphi\|. \qquad _\square
\end{aligned}
$$

\bigskip
\textbf{Step 1.} Denote by $E_{lip}$ the space of all functions
$\varphi \in E$, which are Lipschitz continuous. Here we show that
$P_t\varphi \in C(\mathbb{R})$ for any $\varphi \in E_{lip}$.

Define a sequence $\{\psi_k\}$ of continuous functions
$$
\psi_k(x) = \left\{
\begin{array}{ll}
    1, & \hbox{$x \in [-k,k]$,} \\
    k+1-|x|, & \hbox{$x \in (-k-1,-k) \cup (k,k+1)$,} \\
    0, & \hbox{$x \in (-\infty,-k-1] \cup [k+1,\infty)$.} \\
\end{array}
\right.
$$
Then
$$
|P_t\varphi(x) - P_t\varphi(y)| \leq |O_1| + |O_2| + |O_3|,
$$
where
$$
\begin{aligned}
O_1 &= \mathbb{E}^x \varphi(r_t) \e^{\alpha \int_0^t r_s ds}
\psi_k \left(\int_0^t r_s ds\right) - \mathbb{E}^y \varphi(r_t)
\e^{\alpha
\int_0^t r_s ds} \psi_k \left(\int_0^t r_s ds \right),\\
O_2 &= \mathbb{E}^x \varphi(r_t) \e^{\alpha \int_0^t r_s
ds} \left(1-\psi_k \left( \int_0^t r_s ds \right)\right),\\
O_3 &= \mathbb{E}^y \varphi(r_t) \e^{\alpha \int_0^t r_s ds}
\left(1-\psi_k \left(\int_0^t r_s ds \right)\right).
\end{aligned}
$$

Now we show that $\lim_{y \to x} |O_1| = 0$. To this end write
$$
h_t := \int_0^t r_s ds \qquad {\rm and} \qquad
\zeta_t=(r_t,h_t)^\top
$$
and define a function $\phi_k(\zeta)$ as
$$
\phi_k(r,h) = \varphi(r) \e^{\alpha h} \psi_k(h),
$$
which is Lipschitz continuous with constant $L$, since both
$\phi(r)$ and $\e^{\alpha h} \psi_k(h)$ are Lipschitz continuous.
Denote by $\zeta_t^x$ the value of $\zeta_t$ with initial
condition $\zeta_0 = (x,0)^\top$. We have
$$
|O_1| = |\mathbb{E} [\phi_k(\zeta_t^x) - \phi_k(\zeta_t^y)]| \leq
L \mathbb{E} \| \zeta_t^x - \zeta_t^y \|_2 \leq L \sqrt{\mathbb{E}
\| \zeta_t^x - \zeta_t^y \|_2^2},
$$
where $\| \cdot \|_2$ is the Euclidean norm.

Since
$$
d\zeta_t = \tilde{\mu}(\zeta_t) dt + \tilde{\sigma}(\zeta_t) dW_t
:= \left[
\begin{array}{c}
  a-br_t \\
  r_t \\
\end{array}
\right] dt + \left[
\begin{array}{c}
  \sigma \\
  0 \\
\end{array}
\right] dW_t,
$$
with $\tilde{\mu}$ and $\tilde{\sigma}$ Lipschitz continuous, then
from the mean-square continuity of $\zeta$ (see \cite{Oksendal})
we have
$$
\lim_{y \to x} |O_1| \leq \lim_{y \to x} L \sqrt{\mathbb{E} \|
\zeta_t^x - \zeta_t^y \|_2^2} =0
$$
for all $k \in \mathbb{N}$.

Since we consider $y$ close to $x$, it is now sufficient to show
that $|O_2|$ converges to $0$, as $k \to \infty$, uniformly in $\{
x: |x|<\delta \}$ for any $\delta>0$. We obtain
$$
|O_2| \leq \|\varphi\| \mathbb{E}^x \e^{\frac\alpha b |r_t| +
\alpha h_t} |1-\psi_k(h_t)|,
$$
and from the Schwarz inequality
$$
|O_2| \leq \|\varphi\| \sqrt{\mathbb{E}^x \e^{2\frac\alpha b |r_t|
+ 2\alpha h_t}} \sqrt{\mathbb{P}^x(|h_t| > k)}.
$$
By \eqref{EeNorm} and Chebyshev's inequality,
$$
\begin{aligned}
\mathbb{P}^x(|h_t| > k) &\leq \frac{\mathbb{E}^x |h_t|}{k} \leq
\frac{\mathbb{E}^x \e^{|h_t|}}{k}\\
& \leq \frac{\e^{\frac{\sigma^2 t}{2
b^2}}(1+\e^{\frac{|x|+at}{b}})}{k} \leq \frac{\e^{\frac{\sigma^2
t}{2 b^2}}(1+\e^{\frac{\delta+at}{b}})}{k} \to 0,
\end{aligned}
$$
as $k \to \infty$. In a similar way we show that
$$
\sup_{|x| < \delta} \mathbb{E}^x \e^{2\frac\alpha b |r_t| +
2\alpha h_t} \leq \sup_{|x| < \delta} \sqrt{\mathbb{E}^x
\e^{4\frac\alpha b |r_t|}} \sqrt{\mathbb{E}^x \e^{4\alpha |h_t|}}
< \infty
$$
for every $t \geq 0$. Thus $\lim_{k \to \infty} |O_2| = 0$ and the
convergence is uniform in $\{x: |x|<\delta\}$.

Thus we have
$$
\lim_{y \to x} |P_t\varphi(x) - P_t\varphi(y)| \leq \lim_{k \to
\infty} \lim_{y \to x} (|O_2|+|O_3|) = 0.
$$
Hence, $P_t\varphi$ is continuous for all $t \geq 0$ and $\varphi
\in E_{lip}$.

\bigskip \textbf{Step 2.} We show that $P_t \varphi \in C(\mathbb{R})$ for any $\varphi \in E$.
Let us fix a $\varphi \in E$. As $E_{lip}$ is dense in $E$, there
exists an approximating sequence $\{\varphi_n\}$ such that
$\varphi_n \in E_{lip}$ and $\varphi_n \to \varphi$ in $E$.

Set $\varepsilon > 0$. We have
$$
\begin{aligned}
|P_t & \varphi(x) - P_t\varphi(y)|\\
& \leq |P_t(\varphi-\varphi_n)(x)| + |P_t(\varphi-\varphi_n)(y)| +
|P_t\varphi_n(x) - P_t\varphi_n(y)|\\
& \leq \|P_t(\varphi-\varphi_n)\| \e^{\frac\alpha b |x|} +
\|P_t(\varphi-\varphi_n)\| \e^{\frac\alpha b |y|} +
|P_t\varphi_n(x) - P_t\varphi_n(y)|
\end{aligned}
$$
and from Lemma \ref{l-norm}
$$
\forall \varepsilon >0 \, \exists n_0 \, \forall n>n_0 \quad
\|P_t(\varphi-\varphi_n)\| <\varepsilon.
$$
Furthermore, from Step 1, $P_t \varphi_n \in C(\mathbb{R})$, i.e.
$$
\forall x \in \mathbb{R} \, \exists \delta>0 \, \forall y \in
\mathbb{R} \quad |x-y| < \delta \Rightarrow |P_t\varphi_n(x) -
P_t\varphi_n(y)| <\varepsilon
$$
and therefore
$$
|P_t\varphi(x) - P_t\varphi(y)| <\varepsilon(\e^{\frac{\alpha}{b}
|x|} + \e^{\frac{\alpha}{b} (|x| + \delta)} + 1).
$$

\bigskip \textbf{Step 3.} Here we show that $P_t: E \mapsto E$. For any $\varphi \in E$
write
$$
l(\varphi) = \lim_{|r| \to \infty} |P_t \varphi(r)|
\e^{-\frac\alpha b |r|}.
$$
We need to show that $l(\varphi)=0$. From \eqref{rX} and \eqref{Y}
we have
$$
\begin{aligned}
l(\varphi) \leq & \lim_{|r| \to \infty} \mathbb{E}^r
|\varphi(r_t)| \e^{\alpha\int_0^t r_s ds - \frac \alpha b |r|}\\
= & \lim_{|r| \to \infty} \mathbb{E}^r |\varphi(r_t)|
\e^{\frac{\alpha r}{b} (1-\e^{-bt}) + \frac{\alpha a}{b}(t-
\frac1b (1-\e^{-bt})) + \alpha\sigma Y_t - \frac \alpha b |r|}
\end{aligned}
$$
and from the Schwarz inequality
$$
l(\varphi) \leq \lim_{|r| \to \infty} \e^{\frac{\alpha r}{b}
(1-\e^{-bt}) - \frac \alpha b |r|} \e^{\frac{\alpha a}{b}(t-
\frac1b (1-\e^{-bt}))} \sqrt{\mathbb{E} \e^{2\alpha\sigma Y_t}}
\sqrt{\mathbb{E}^r \varphi^2(r_t)}.
$$
Hence
$$
l^2(\varphi) \leq g^2(t) \lim_{|r| \to \infty} \e^{2\frac{\alpha
r}{b} (1-\e^{-bt}) - 2\frac \alpha b |r|} \mathbb{E}^r
\varphi^2(r_t),
$$
where
$$
g(t) := \e^{\frac{\alpha a}{b}(t- \frac1b (1-\e^{-bt}))}
\sqrt{\mathbb{E} \e^{2\alpha\sigma Y_t}}.
$$
Clearly $g(t)<\infty$ for any fixed $t \geq 0$, by \eqref{Y-dist}
and \eqref{EeNorm}.

Set $\varepsilon>0$. Since $\varphi \in E$, then there exists a
$\delta>0$ such that
$$
|\varphi(r_t)| \leq \varepsilon \e^{\frac \alpha b |r_t|}
$$
in a set $\{ \omega: |r_t| \geq \delta \}$ and therefore
$$
\mathbb{E}^r \varphi^2(r_t) \leq \varepsilon^2 \mathbb{E}^r
\e^{\frac{2\alpha}{b} |r_t|} + \|\varphi\|_{\delta}^2,
$$
where $\|\cdot\|_{\delta}$ is the supremum norm over
$\{|x|<\delta\}$. Clearly $\|\varphi\|_{\delta}<\infty$ for every
$\varphi \in E$. From \eqref{rX},
$$
\begin{aligned}
l^2(\varphi) \leq & \varepsilon^2 g^2(t) \e^{\frac{2\alpha
a}{b^2}(1-\e^{-bt})} \mathbb{E} \e^{\frac{2\alpha \sigma}{b}|X_t|}
\lim_{|r| \to \infty}
\e^{2\frac{\alpha}{b} (1-\e^{-bt})(r - |r|)}\\
&+ \|\varphi\|_{\delta}^2 g^2(t) \lim_{|r| \to \infty}
\e^{2\frac{\alpha}{b} (r(1-\e^{-bt}) - |r|)}
\end{aligned}
$$
and
$$
l^2(\varphi) \leq  2 \varepsilon^2 g^2(t) \e^{\frac{2\alpha
a}{b^2}} \e^{\frac{2\alpha^2 \sigma^2}{b^3}(1-\e^{-2bt})} \leq 2
\varepsilon^2 g^2(t) \e^{\frac{2\alpha a}{b^2} + \frac{2\alpha^2
\sigma^2}{b^3}}
$$
by \eqref{X-dist} and \eqref{EeNorm}. Hence, $P_t \varphi \in E$.

\bigskip \textbf{Step 4.} Clearly $P_0=I$ and $P_t P_s =
P_{t+s}$ holds since $(r_t)$ is a Markov process. We need to show
strong continuity of $(P_t)$, i.e
\begin{equation}\label{strong-conv}
\lim_{t \downarrow 0} \| P_t \varphi - \varphi \| = 0, \qquad
\forall \varphi \in E.
\end{equation}

Taking into account Lemma \ref{l-norm} and the Banach--Steinhaus
theorem, it is enough to show \eqref{strong-conv} for $\varphi \in
C_0(\mathbb{R})$. For such a $\varphi$ we have
$$
\| P_t\varphi - \varphi \| = \sup_{r \in \mathbb{R}} |\mathbb{E}^r
\varphi(r_t) \e^{\alpha \int_0^t r_s ds} - \varphi(r)| \e^{-\frac
\alpha b |r|} \leq O_1 + O_2 + O_3,
$$
where
$$
\begin{aligned}
O_1 &= \sup_{r \in \mathbb{R}} |\mathbb{E}^r (\e^{\alpha \int_0^t
r_s ds} - 1) \varphi(r)| \e^{-\frac \alpha b |r|},\\
O_2 &= \sup_{r \in \mathbb{R}} |\mathbb{E}^r
\varphi(r_t) - \varphi(r)| \e^{-\frac \alpha b |r|},\\
O_3 &= \sup_{r \in \mathbb{R}} |\mathbb{E}^r (\e^{\alpha \int_0^t
r_s ds}-1)(\varphi(r_t) - \varphi(r))| \e^{-\frac \alpha b |r|}.
\end{aligned}
$$
Since $\varphi \in C_0(\mathbb{R})$, then there exists a
$\delta>0$ such that $supp \, \varphi \subset (-\delta,\delta)$
and
$$
\lim_{t \downarrow 0} O_1 \leq \|\varphi\|_{\delta} \lim_{t
\downarrow 0} \sup_{|r| < \delta} |\e^{f(r,t)}-1| \e^{-\frac
\alpha b |r|}
$$
where
$$
\begin{aligned}
f(r,t) = & \frac{\alpha r}{b} (1-\e^{-bt}) + \frac{\alpha a}{b}
\left( t- \frac1b (1-\e^{-bt}) \right)\\
& + \frac{\alpha^2\sigma^2}{2b^2} \left( t-\frac{3}{2b} + \frac2b
\e^{-bt} - \frac{1}{2b} \e^{-2bt} \right).
\end{aligned}
$$
Since $\lim_{t \to 0}f(r,t) = 0$ uniformly in $\{|r|<\delta \}$,
we easily verify that $\lim_{t \downarrow 0} O_1 = 0$.

Set $\varepsilon>0$. Recall, that every continuous function on a
compact set is uniformly continuous. Thus there exists a $\rho
>0$ such that
$$
\begin{aligned}
O_2 &\leq \varepsilon \sup_{r \in \mathbb{R}}
\mathbb{P}(|r_t-r|<\rho) \e^{-\frac \alpha b |r|} + \sup_{r \in
\mathbb{R}} \int_{\{ |r_t-r|\geq\rho \}}
|\varphi(r_t)-\varphi(r)| \e^{-\frac \alpha b |r|} d\mathbb{P}\\
&\leq \varepsilon + 2 \|\varphi\|_\infty \sup_{r \in \mathbb{R}}
\mathbb{P}(|r_t-r|\geq\rho) \e^{-\frac \alpha b |r|},
\end{aligned}
$$
where $\|\cdot\|_\infty$ is the supremum norm. Hence
$$
\begin{aligned}
\lim_{t \downarrow 0} O_2 &\leq \varepsilon + 2 \|\varphi\|_\infty
\lim_{t \downarrow 0} \sup_{r \in \mathbb{R}}
\mathbb{P}(|r_t-r|\geq\rho) \e^{-\frac \alpha b |r|}\\
&=\varepsilon + 2 \|\varphi\|_\infty \lim_{t \downarrow 0} \sup_{r
\in \mathbb{R}} \e^{-\frac \alpha b |r|} (1 + N_{0,1}(r_1) -
N_{0,1}(r_2))
\end{aligned}
$$
where $N_{0,1}(\cdot)$ is a normal distribution function of
$\mathcal{N}(0,1)$ and
$$
r_1 = \sqrt{2b}\tfrac{(r-\frac
ab)(1-\e^{-bt})-\rho}{\sigma\sqrt{1-\e^{-2bt}}}, \qquad r_2 =
\sqrt{2b}\tfrac{(r-\frac
ab)(1-\e^{-bt})+\rho}{\sigma\sqrt{1-\e^{-2bt}}}.
$$
The supremum is attained at $r = \pm \infty$ or $r=\hat{r}(t)
<\infty$, but with a possible infinite limit, i.e. $\lim_{t
\downarrow 0} |\hat{r}(t)| \leq \infty$. In all this cases we
obtain
$$
\lim_{t \downarrow 0} \sup_{r \in \mathbb{R}} \e^{-\frac \alpha b
|r|} (1 + N_{0,1}(r_1) - N_{0,1}(r_2)) =0.
$$
Hence, taking $\varepsilon \to 0$, $\lim_{t \downarrow 0} O_2 =0$.

Finally, from the Schwartz inequality
$$
\lim_{t \downarrow 0} O_3 \leq 2 \|\varphi\|_\infty \lim_{t
\downarrow 0} \sup_{r \in \mathbb{R}} \sqrt{\e^{-2 \frac \alpha b
|r|} \mathbb{E}^r (\e^{\alpha\int_0^t r_s ds}-1)^2},
$$
where we can easily derive an analytic formula of $\mathbb{E}^r
(\e^{\alpha\int_0^t r_s ds}-1)^2$. Then we take into consideration
all the possible realization of supremum $\hat{r}$, as above, and
we get $\lim_{t \downarrow 0} O_3 =0$.

\section*{Appendix $B$ - proof of $C_0$-semigroup property of $(P_t)$ in case of invariant interval model}

In order to prove that $P_t \colon E \mapsto E$ we need to show
that $P_t \varphi \in UC((a,b))$ and this may be done in a similar
way to the proof of continuity of $(P_t)$ in Appendix $A$. Since
$P_0=I$ and $P_t P_s = P_{t+s}$ clearly hold, so in order to prove
$C_0$-semigroup property of $(P_t)$ we need to show that
\eqref{strong-conv} holds for $E=UC((a,b))$ equipped with the
supremum norm
$$
\|\varphi\| = \sup_{r \in (a,b)} |\varphi(r)|.
$$

Let $\varphi \in E$, then
$$
\| P_t\varphi - \varphi \| = \sup_{r \in (a,b)} |\mathbb{E}^r
\varphi(r_t) \e^{\alpha \int_0^t r_s ds} - \varphi(r)| \leq
O_1+O_2+O_3,
$$
where
$$
\begin{aligned}
O_1 &= \sup_{r \in (a,b)} |\mathbb{E}^r (\e^{\alpha \int_0^t r_s
ds} - 1) \varphi(r)|,\\
O_2 &= \sup_{r \in (a,b)} |\mathbb{E}^r \varphi(r_t) -
\varphi(r)|,\\
O_3 &= \sup_{r \in (a,b)} |\mathbb{E}^r (\e^{\alpha \int_0^t r_s
ds}-1)(\varphi(r_t) - \varphi(r))|.
\end{aligned}
$$
Since $\|\varphi\|<\infty$ for every $\varphi \in E$, then
$$
\begin{aligned}
\lim_{t \downarrow 0} O_1 & \leq \|\varphi\| \lim_{t \downarrow 0}
\mathbb{E}^r |\e^{\alpha \int_0^t r_s ds} - 1|\\
& \leq \|\varphi\| \lim_{t \downarrow 0} ( \max\{|\e^{\alpha at} -
1|, |\e^{\alpha bt} - 1| \}) =0
\end{aligned}
$$
and
$$
\lim_{t \downarrow 0} O_3 \leq 2 \|\varphi\| \lim_{t \downarrow 0}
( \max \{ |\e^{\alpha at} - 1|, |\e^{\alpha bt} - 1| \}) =0.
$$

Let $E_{lip}$ be the space of all Lipschitz continuous functions
$\varphi \in E$. Let $\varphi \in E_{lip}$ and let $L$ be the
Lipschitz constant of $\varphi$. Then we have
$$
\begin{aligned}
\| P_t \varphi &- \varphi \| = \sup_{r \in (a,b)} |\mathbb{E}^r
\varphi(r_t) - \varphi(r)|\\
&\leq L \sup_{r \in (a,b)} \mathbb{E}^r |r_t - r| = L \sup_{r \in
(a,b)} \mathbb{E}^r \left| \int_0^t
\mu(r_s)ds + \int_0^t \sigma(r_s)dW_s \right| \\
&\leq L \sup_{r \in (a,b)} \sqrt{\mathbb{E}^r \left( \int_0^t
\mu(r_s)ds + \int_0^t \sigma(r_s)dW_s \right)^2}\\
&\leq L \sup_{r \in (a,b)} \sqrt{2 \mathbb{E}^r \left( \int_0^t
\mu(r_s)ds \right)^2 + 2 \mathbb{E}^r \left( \int_0^t
\sigma(r_s)dW_s \right)^2}\\
&\leq L \sqrt{2 t^2 \sup_{r \in (a,b)} |\mu(r)|^2  + 2 t \sup_{r
\in (a,b)} |\sigma(r)|^2}.
\end{aligned}
$$
Hence, $\lim_{t \downarrow 0} \| P_t \varphi - \varphi \| = 0$ for
every $\varphi \in E_{lip}$. Since $E_{lip}$ is dense in $E$, we
conclude by the Banach--Steinhaus theorem.

\end{document}